\documentclass[11pt]{article}

\usepackage{amsmath, amstext, amsgen, amsbsy, amsopn, amsfonts, amssymb, graphicx,psfrag,authblk}
\usepackage[normalem]{ulem}
\usepackage{pdfsync}
\usepackage{epstopdf}
\usepackage{color}
\usepackage[thinlines]{easytable}

\def\Z{\mathbb{Z}}

\newtheorem{theorem}{Theorem}
\newtheorem{proposition}[theorem]{Proposition}
\newtheorem{lemma}[theorem]{Lemma}
\newtheorem{corollary}[theorem]{Corollary}

\begin{document}

\markboth{Crans, Hoste, Mellor and Shanahan}{Finite $n$-quandles of torus and two-bridge links}

\title{Finite $n$-quandles of torus and two-bridge links}

\author{Alissa S. Crans, Blake Mellor, and Patrick D. Shanahan}
\affil{{\em Loyola Marymount University} \\ {\em 1 LMU Drive} \\ {\em Los Angeles, CA 90045}}
\author{Jim Hoste}
\affil{{\em Pitzer College} \\ {\em 1050 N. Mills Avenue} \\ {\em Claremont, CA 91711}}

\date{}
\maketitle
\begin{abstract} We compute Cayley graphs and automorphism groups for all finite $n$-quandles of two-bridge and torus knots and links, as well as torus links with an axis.  \end{abstract}

\section{Introduction}
Associated to every oriented knot and link $L$ is its  fundamental quandle $Q(L)$. Except for the unknot and Hopf link, these quandles are infinite, but for each integer $n>1$  a certain quotient of $Q(L)$, called the  $n$-quandle of $L,$ and denoted by $Q_n(L)$, may be finite. From results of Joyce \cite{JO2, JO} and Winker \cite{WI}, if the $n$-quandle $Q_n(L)$ is finite, then $\widetilde M_n(L)$, the $n$-fold cyclic branched cover of $S^3$ branched over $L$, has finite fundamental group. It was conjectured by Przytycki, and recently proven by Hoste and Shanahan \cite{HS2}, that the converse is also true.  Using this result, together with Dunbar's \cite{DU} classification of all geometric, non-hyperbolic 3-orbifolds, a complete list of all  knots and links in $S^3$ with finite $n$-quandle for some $n$ was given in \cite{HS2}.  The links are listed in Table \ref{linktable} (reproduced from \cite{HS2}) as they are given by Dunbar and include all two-bridge links, some torus  links, some torus links with ``axis,'' and some Montesinos links.
This paper is the second in a series of  papers, beginning with \cite{HS1},  to give detailed descriptions of these finite $n$-quandles with the ultimate goal being a tabulation of all finite quandles that appear as $n$-quandles of links for some $n$.  Here we describe the  Cayley graphs and automorphism groups of the finite $n$-quandles associated to  the two-bridge links, torus links, and torus links with axis. These links correspond to the links in the first three rows and the first entry in the fourth row of Table \ref{linktable}. (We warn the reader that it is not obvious that the first entry in the fourth row of this table represents all possible 2-bridge links.) The Montesinos links considered in \cite{HS1} appear as the third entry in the fourth row.  The four remaining links in Table \ref{linktable} are infinite families which, like the Montesinos links, require substantially more analysis than the cases considered here. We intend to consider these families in future work.

\begin{table}[htbp]
$$
\begin{array}{ccc}
\includegraphics[width=1.25in,trim=0 0 0 0,clip]{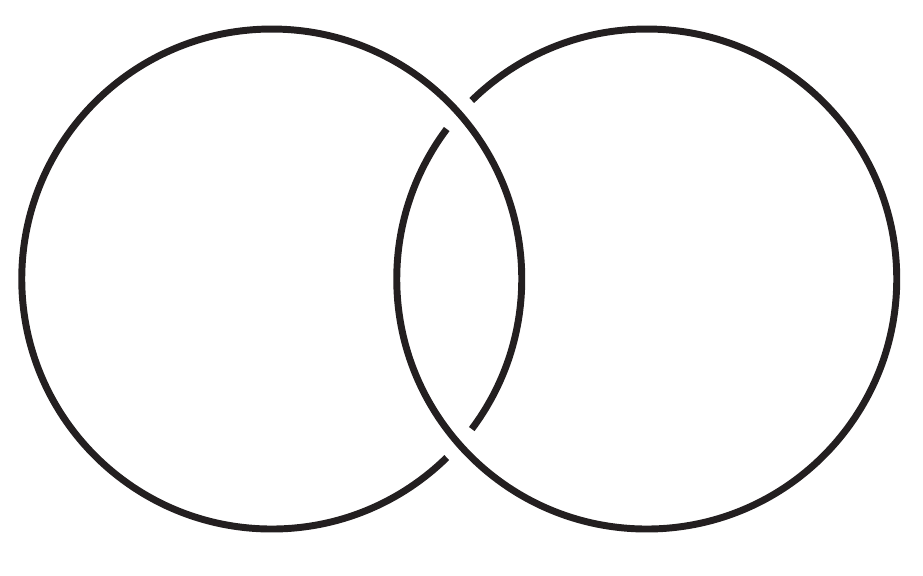}   & \includegraphics[width=1.25in,trim=0 0 0 0,clip]{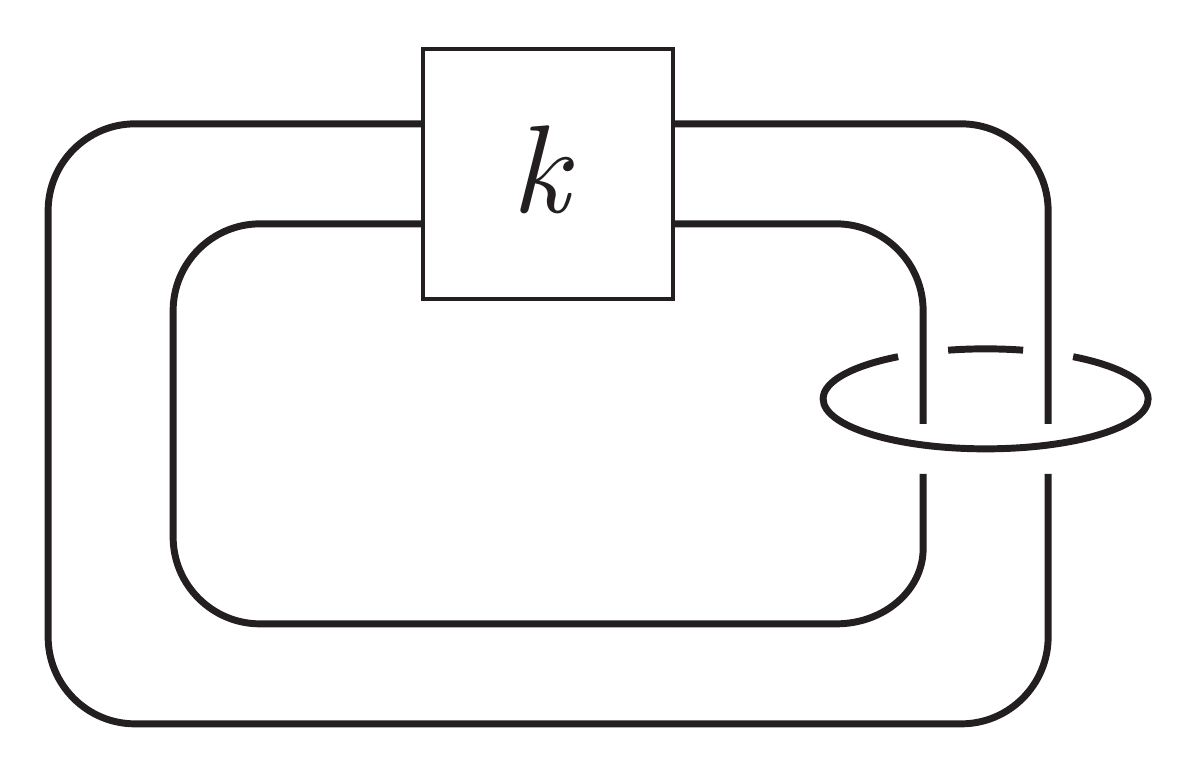}  & \includegraphics[width=1.0in,trim=0 0 0 0,clip]{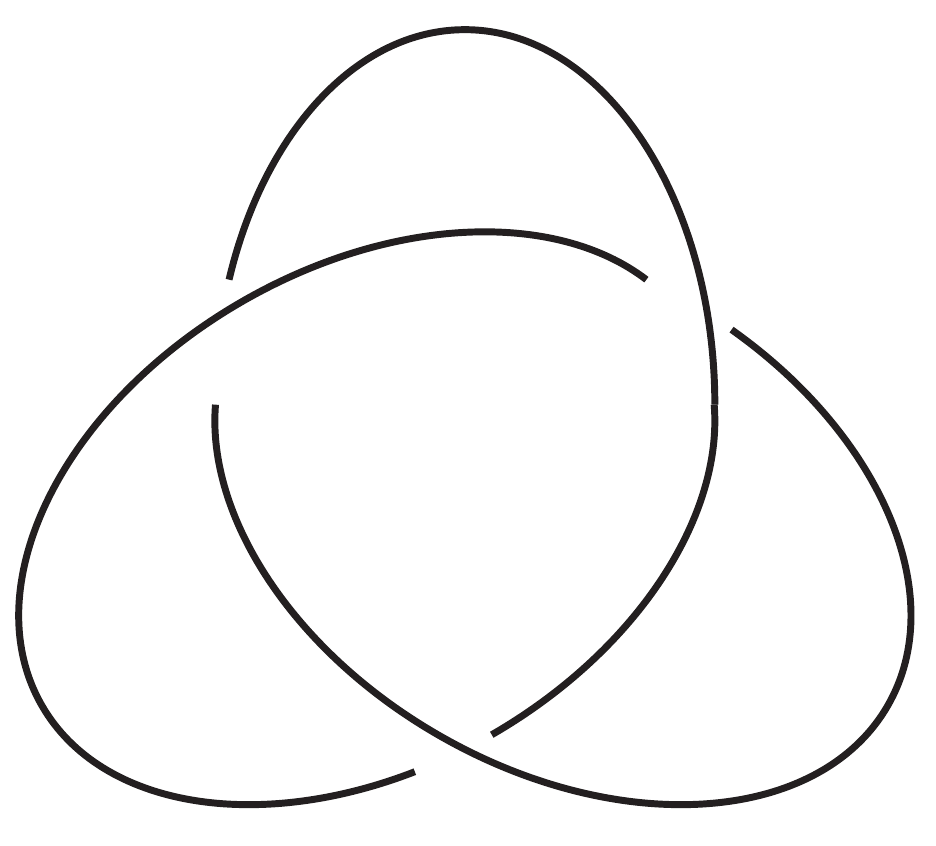} \\
\scriptstyle n > 1 & \scriptstyle k \neq 0,\  n=2& \scriptstyle n=3, 4, 5 \\
\\
\includegraphics[width=1.0in,trim=0 0 0 0,clip]{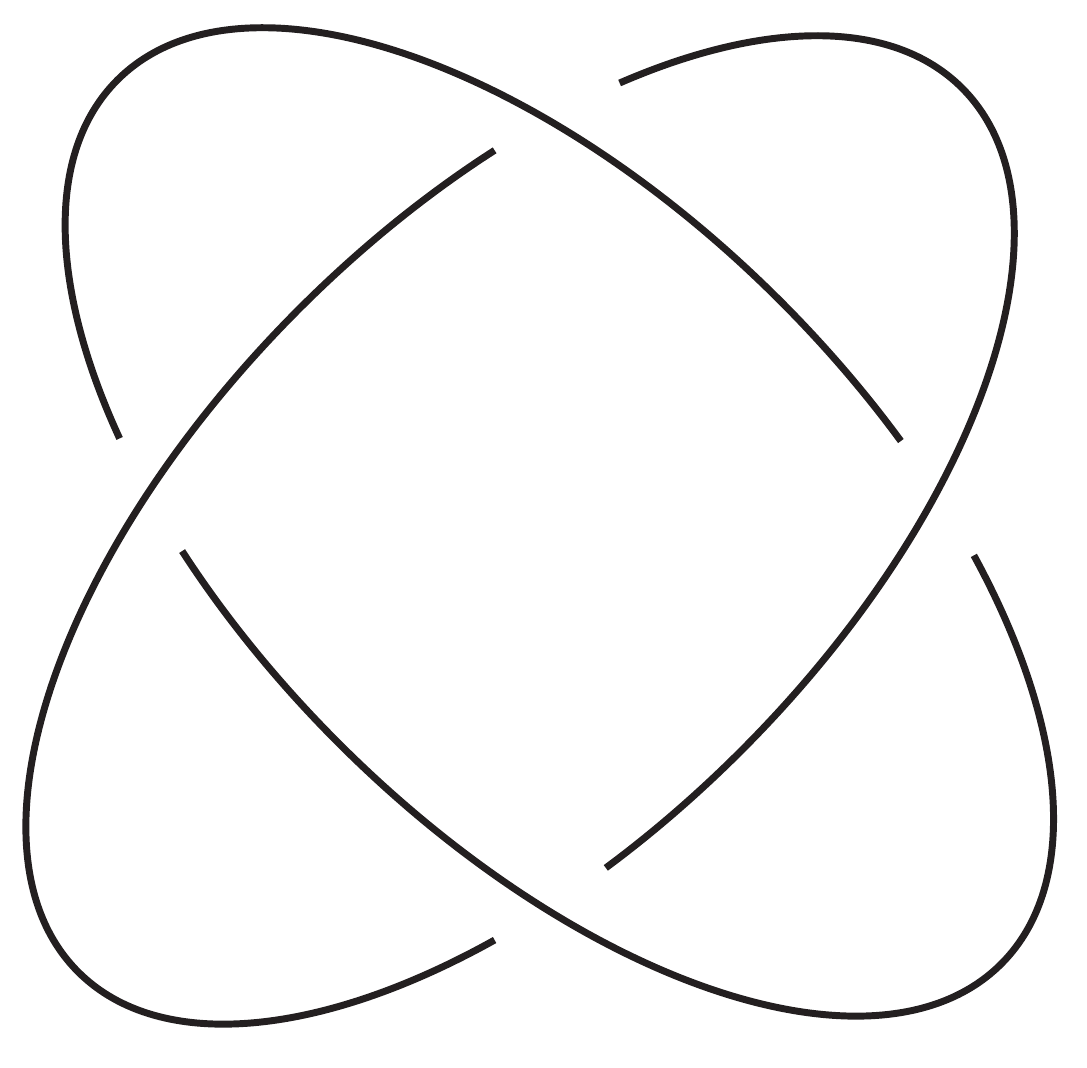}  & \includegraphics[width=1.0in,trim=0 0 0 0,clip]{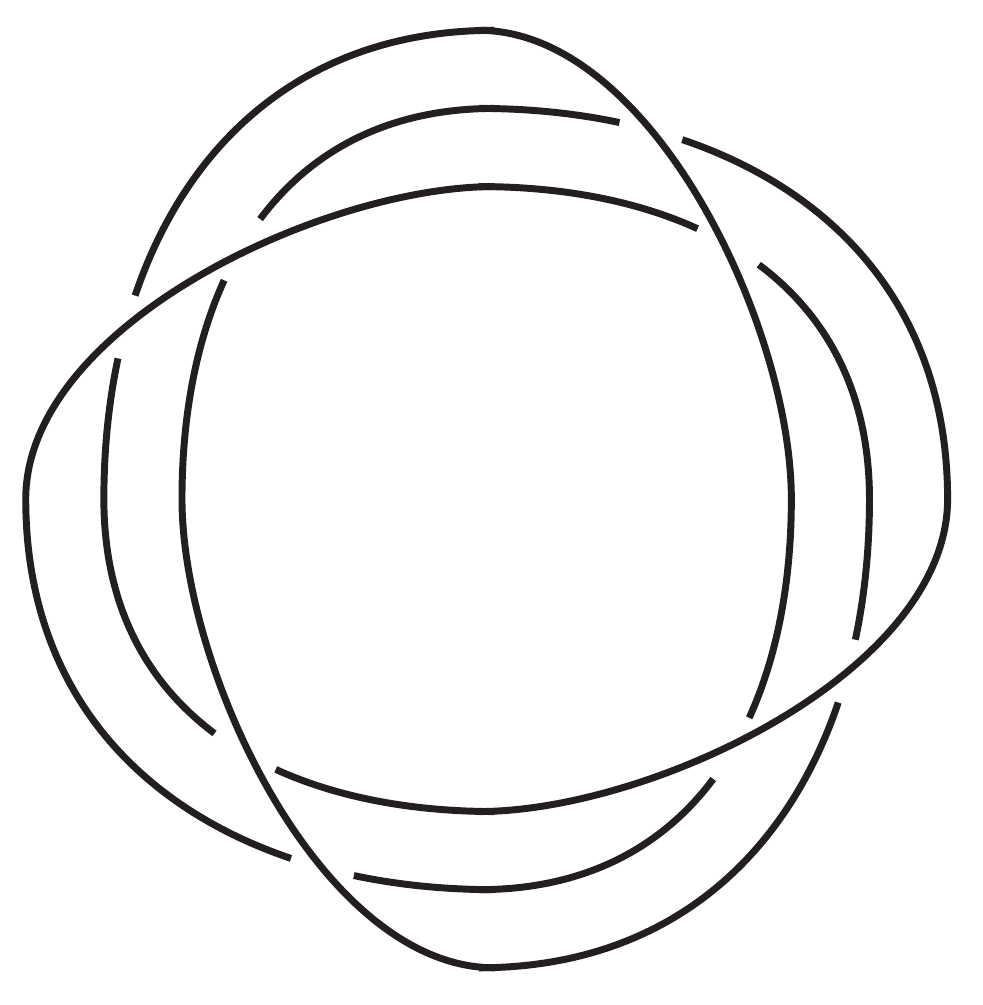}  & \includegraphics[width=1.0in,trim=0 0 0 0,clip]{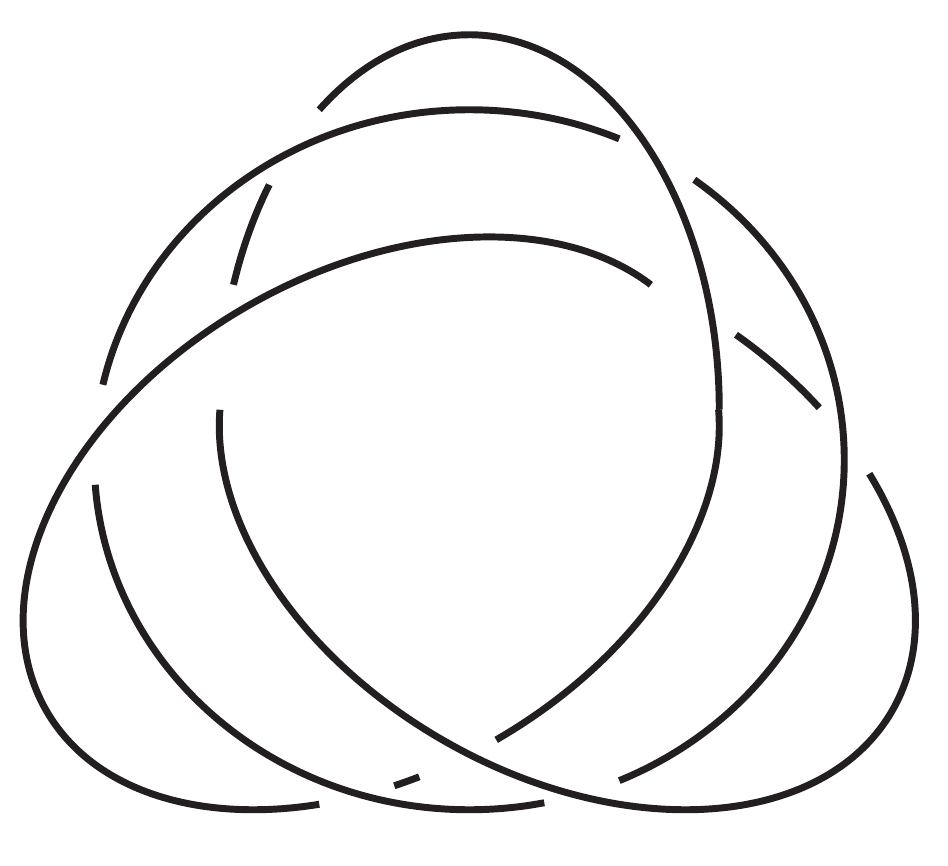} \\
\scriptstyle n =3 & \scriptstyle n=2& \scriptstyle n=2 \\
\\
\includegraphics[width=1.0in,trim=0 0 0 0,clip]{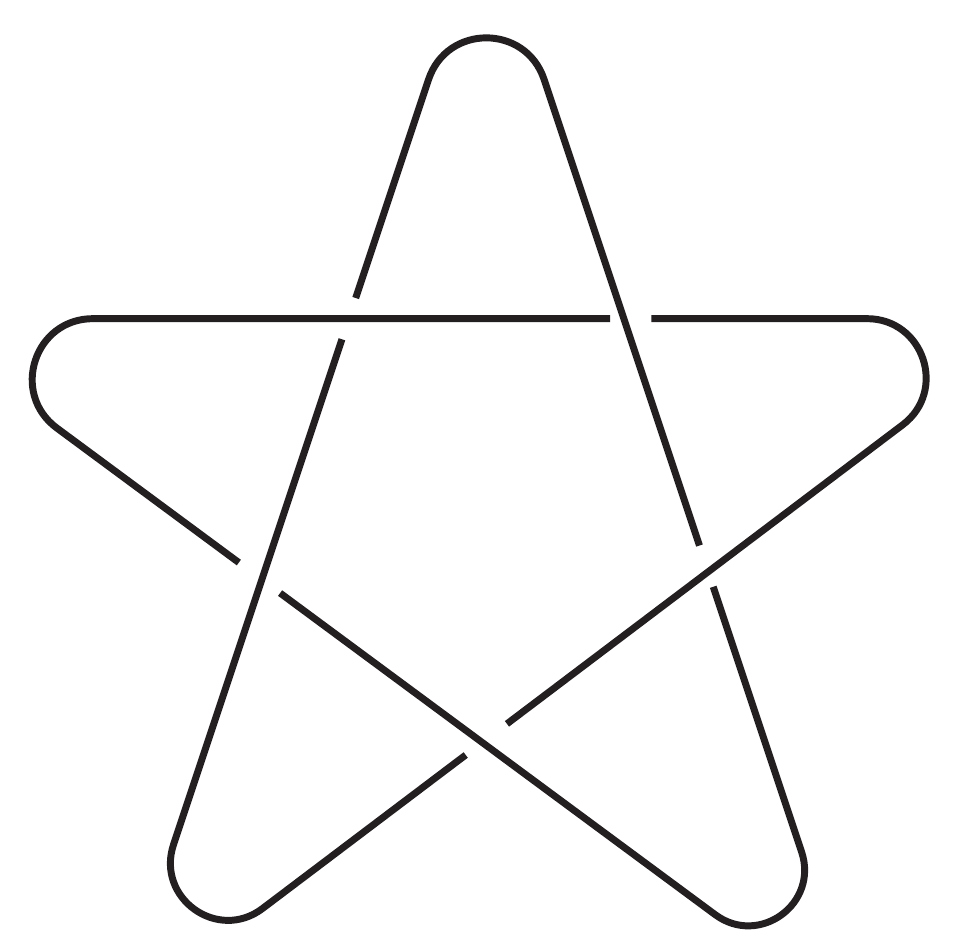}  & \includegraphics[width=1.0in,trim=0 0 0 0,clip]{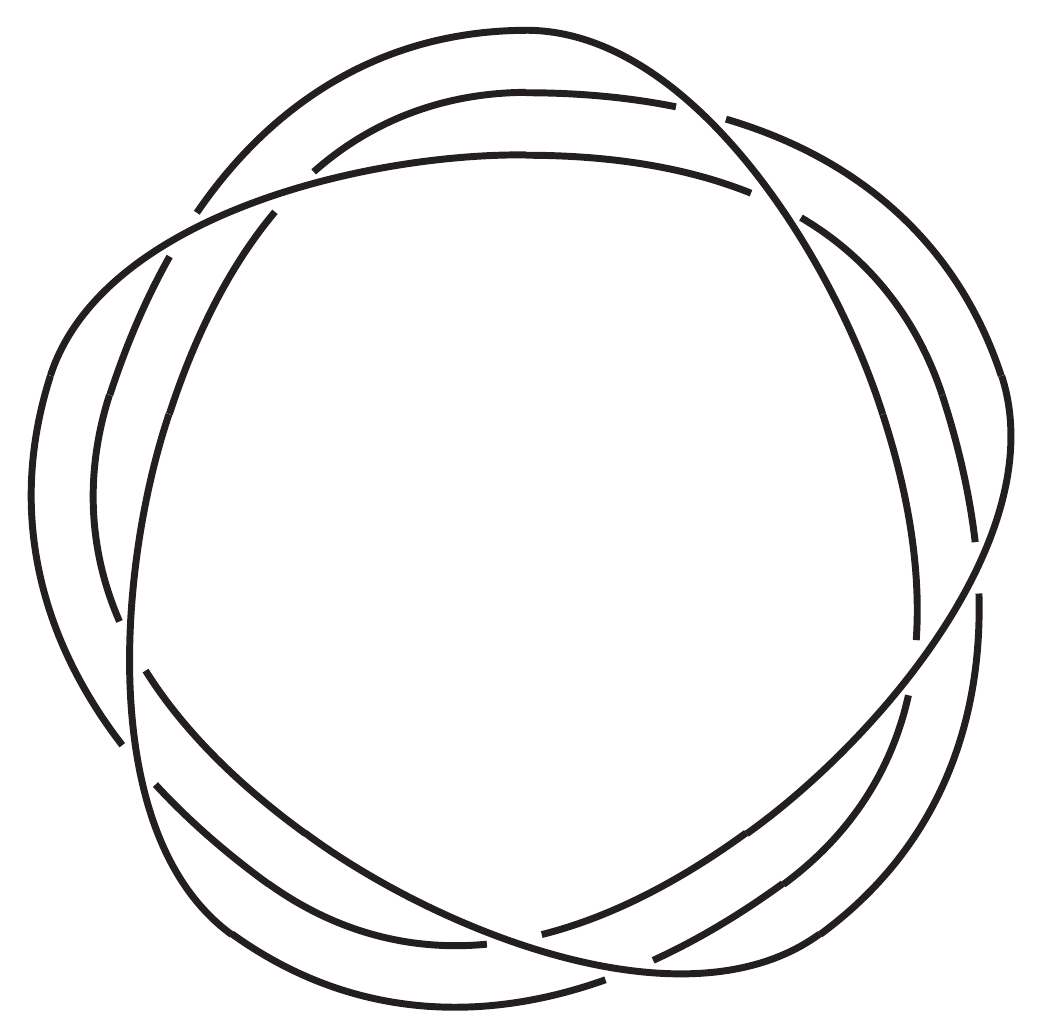}  & \includegraphics[width=1.25in,trim=0 0 0 0,clip]{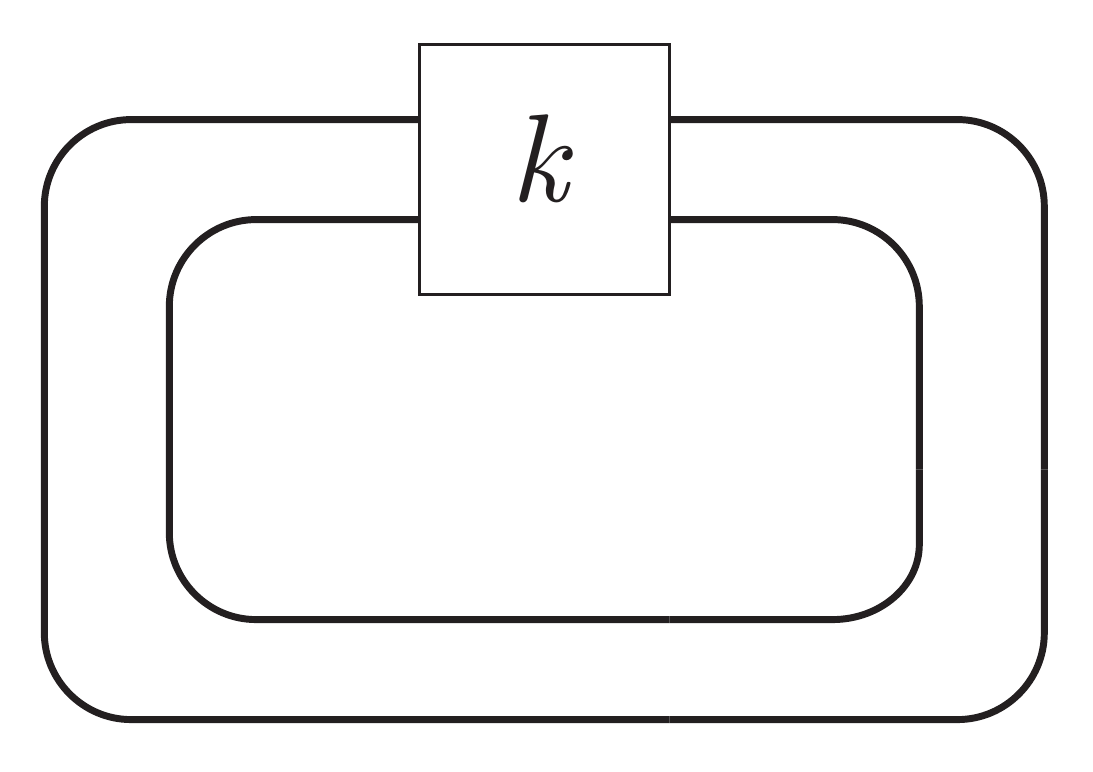} \\
\scriptstyle n =3 & \scriptstyle n=2& \scriptstyle k\neq0,\ n=2 \\
\\
\includegraphics[width=1.25in,trim=0 0 0 0,clip]{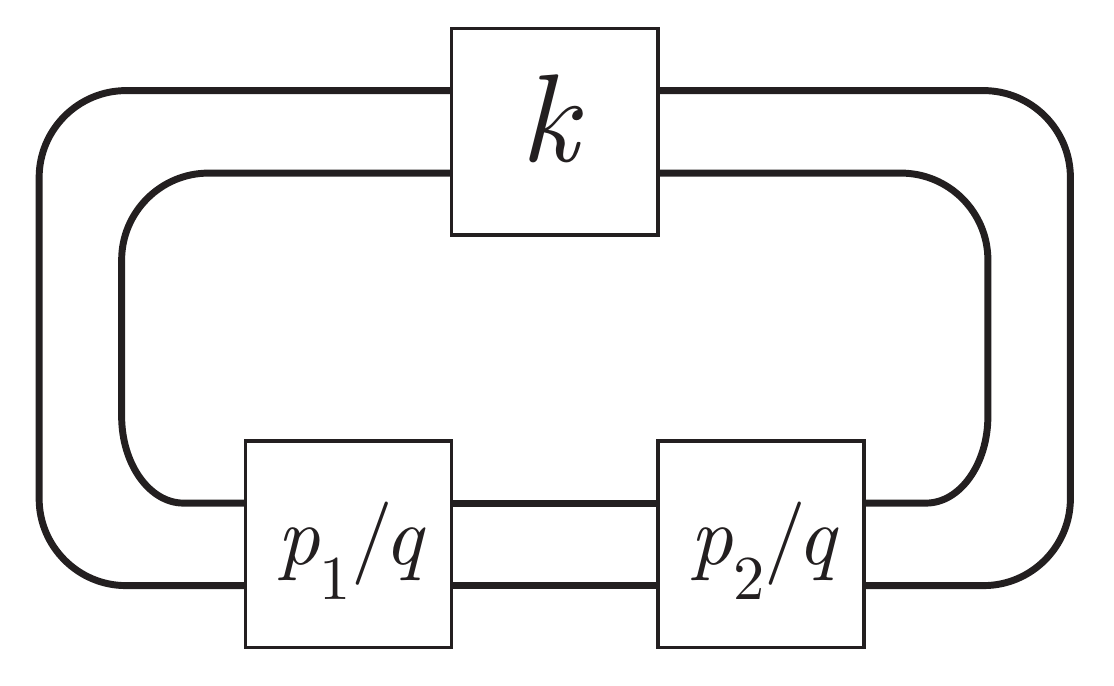}  & \includegraphics[width=1.15in,trim=0 5pt 0 0,clip]{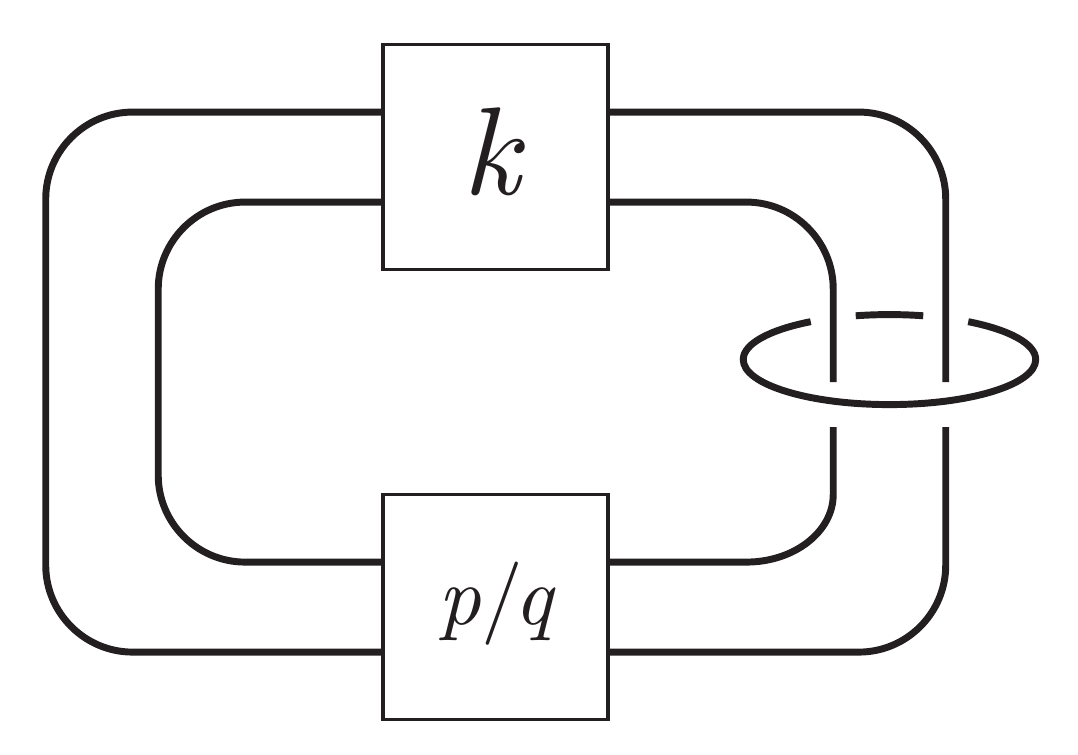} &  \includegraphics[width=1.65in,trim=0 0 0 0,clip]{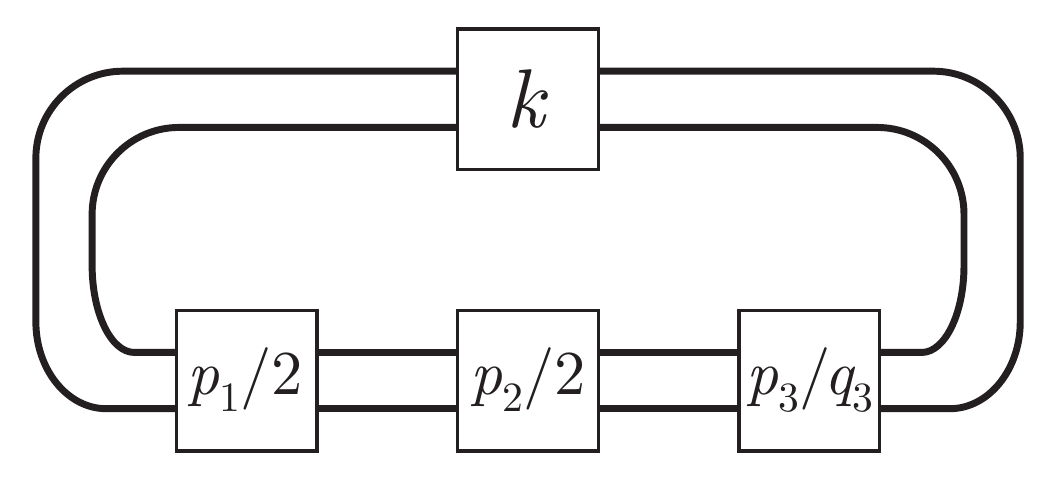} \\
\scriptstyle k+p_1/q+p_2/q \neq 0,\ n =2  &\scriptstyle n=2& \scriptstyle k+p_1/2+p_2/2+p_3/q_3 \neq 0,\ n =2 \\
\\
\includegraphics[width=1.65in,trim=0pt 0pt 0pt 0pt,clip]{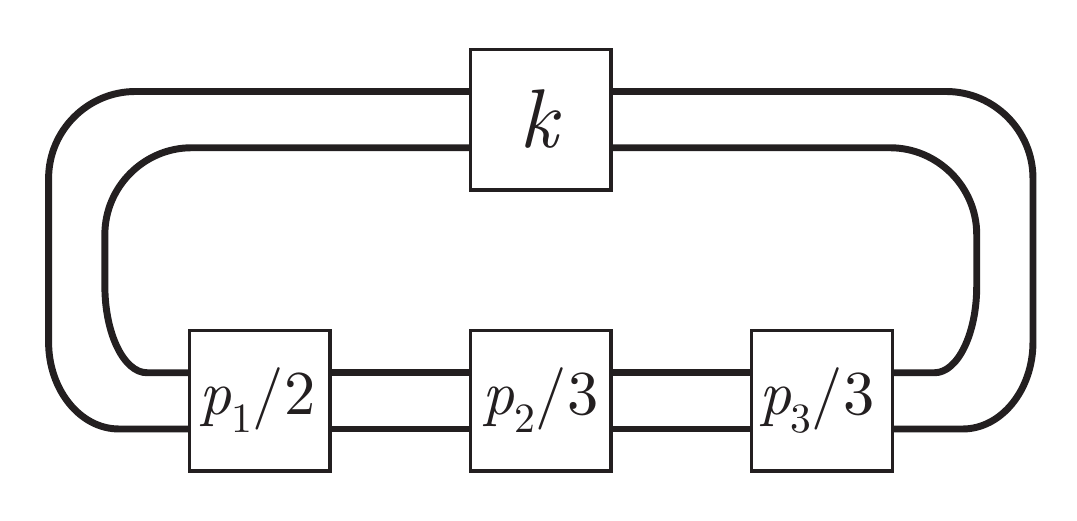}  & \includegraphics[width=1.65in,trim=0pt 0pt 0pt 0pt,clip]{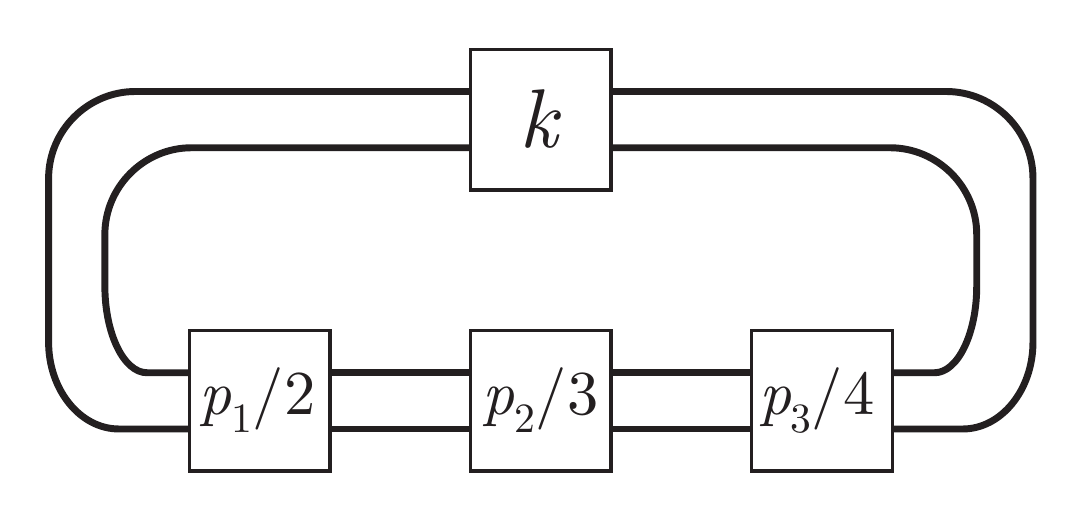}  & \includegraphics[width=1.65in,trim=0pt 0pt 0pt 0pt,clip]{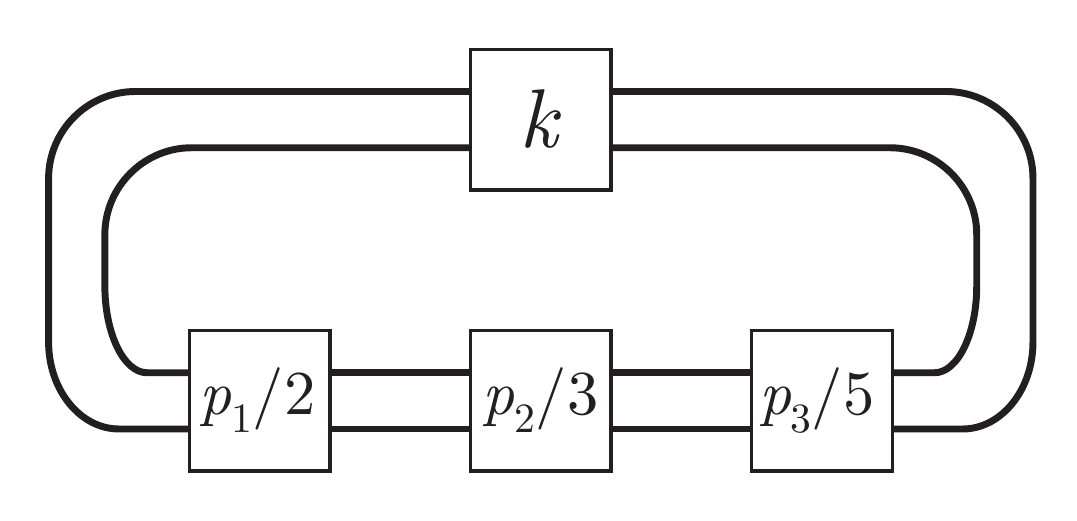}\\
\scriptstyle k+p_1/2+p_2/3+p_3/3 \neq 0,\ n =2  & \scriptstyle k+p_1/2+p_2/3+p_3/4 \neq 0,\ n =2  & \scriptstyle k+p_1/2+p_2/3+p_3/5 \neq 0,\ n =2
\end{array}
$$
\caption{Links $L \in \mathbb{S}^3$ with finite $Q_n(L)$}
\label{linktable}
\end{table}

\section{Link quandles}

We begin with a review of the definition of the fundamental quandle of a link and its associated $n$-quandles. We refer the reader to \cite{FR}, \cite{JO2}, \cite{JO}, and \cite{WI} for more detailed information.

A {\it  quandle} is a set $Q$ equipped with two binary operations $\rhd$ and $\rhd^{-1}$ that satisfy the following three axioms:
\begin{itemize}
\item[\bf A1.] $x \rhd x =x$ for all $x \in Q$.
\item[\bf A2.] $(x \rhd y) \rhd^{-1} y = x = (x \rhd^{-1} y) \rhd y$ for all $x, y \in Q$.
\item[\bf A3.] $(x \rhd y) \rhd z = (x \rhd z) \rhd (y \rhd z)$ for all $x,y,z \in Q$.
\end{itemize}

Each element $x\in Q$ defines a map $S_x:Q \to Q$ by $S_x(y)=y \rhd x$. The axiom A2 implies that each $S_x$ is a bijection and the axiom A3 implies that each $S_x$ is a quandle homomorphism, and therefore an automorphism. We call $S_x$ the {\it point symmetry at $x$}.

It is important to note that the operation $\rhd$ is, in general, not associative. In order to clarify the ambiguity caused by lack of associativity, we adopt the exponential notation introduced by  Fenn and Rourke in \cite{FR} and denote $x \rhd y$ as $x^y$ and $x \rhd^{-1} y$ as $x^{\bar y}$. With this notation, $x^{yz}$ will be taken to mean $(x^y)^z=(x \rhd y)\rhd z$ whereas $x^{y^z}$ will mean $x\rhd (y \rhd z)$. 

The following lemma from \cite{FR}, which describes how to re-associate a product in an $n$-quandle given by a presentation, will be used repeatedly in this paper. 
\begin{lemma} If $a^u$ and $b^v$ are elements of a quandle, then
$$\left(a^u \right)^{\left(b^v \right)}=a^{u \bar v b v} \ \ \ \ \mbox{and}\ \ \ \ \left(a^u \right)^{\overline{\left(b^v \right)}}=a^{u \bar v \bar b v}.$$
\label{leftassoc}
\end{lemma}
Using Lemma~\ref{leftassoc}, elements in a quandle given by a presentation $\langle S \mid R \rangle$ can be represented as equivalence classes of expressions of the form $a^w$ where $a$ is a generator in $S$ and $w$ is a word in the free group on $S$ (with $\bar x$ representing the inverse of $x$).

If $n$ is a natural number, a quandle $Q$ is an {\em $n$-quandle} if  $x^{y^n} =x$ for all $x,y \in Q$, where by $y^n$ we mean  $y$ repeated $n$ times. In other words, each point symmetry $S_x$ has order dividing $n$. A {\em trivial} quandle  is one where $x^y=x^{\bar y}=x$ for all $x, y\in Q$ or, equivalently,  a trivial quandle is a 1-quandle. 
A $2$-quandle is also called an {\it involutory} quandle. Note that a  quandle is involutory  if and only if $\rhd = \rhd^{-1}$, that is, each point symmetry is an involution.

\begin{figure}[h]
$$\includegraphics[height=1.25in]{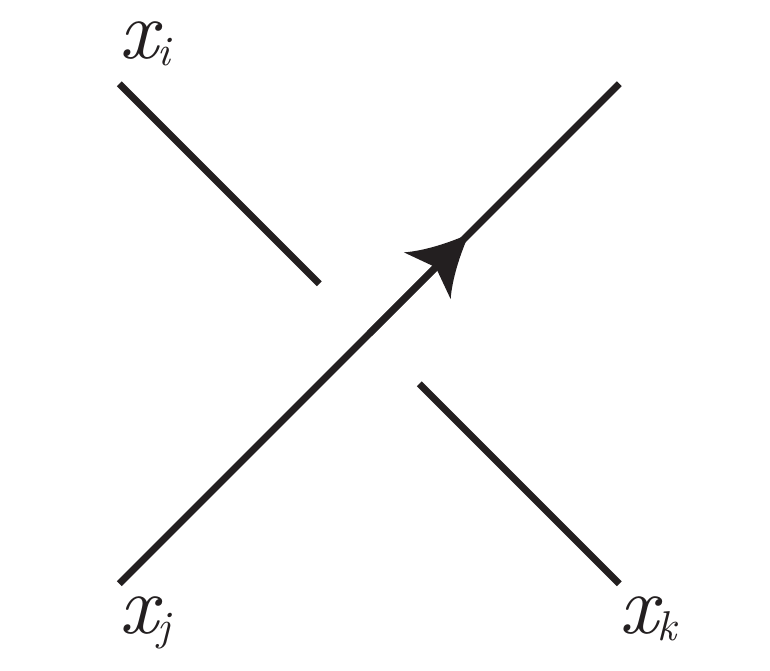}$$
\caption{The relation $x_i=x_k^{x_j}$ at a crossing.}
\label{crossing}
\end{figure}

If $L$ is an oriented knot or link in $S^3$, then a presentation of its fundamental quandle, $Q(L)$, can be derived from a regular diagram $D$ of $L$. See Joyce~\cite{JO}. This process mimics the Wirtinger algorithm. Namely, assign a quandle generator $x_1, x_2, \dots , x_n$ to each arc of $D$, then at each crossing introduce the relation $x_i=x_k^{ x_j}$ as shown in Figure~\ref{crossing}. It is easy to check that the three Reidemeister moves do not change the quandle given by this presentation so that the quandle is indeed an invariant of the oriented link.

The fundamental quandle of a link depends on the choice of orientation. If $L$ is an oriented link, let $-L$ be its {\em reverse} obtained by reversing the orientation of all of its components and let $L^*$ be its {\em obverse}, the mirror image of $L$. The {\em inverse} of $L$ is defined to be $-L^*$. If $D$ is a diagram of an oriented link and $D'$ is obtained from $D$ by both reflecting through a plane perpendicular to the plane of projection and reversing all orientations, then it is not hard to see that both $D$ and $D'$ produce exactly the same presentation for the fundamental quandle. Thus, the fundamental quandle is the same for an oriented link $L$ and its inverse $-L^*$. 

Given an oriented link $L$ and a presentation $\langle S \,|\, R\rangle$ of $Q(L)$, a presentation of  the quotient $n$-quandle $Q_n(L)$ is obtained by adding the relations $x^{y^n}=x$ for every pair of distinct generators $x$ and $y$. As with the fundamental quandle, $Q_n(L)$ depends on the choice of orientation but $Q_n(L) = Q_n(-L^*)$. In the case $n=2$, the relation $x_i=x_k^{x_j}$ is equivalent to the relation $x_k=x_i^{x_j}$. Hence, the 2-quandle does not depend on the orientation of the link, so $Q_2(L) = Q_2(-L) = Q_2(L^*) = Q_2(-L^*)$. In general, one should expect a significant loss of information in passing from $Q(L)$ to $Q_n(L)$. 

The following result is a simple corollary of the main theorem of \cite{HS2}.

\begin{proposition} The fundamental quandle of an oriented link $L$ is finite if and only if $L$ is either the unknot or the Hopf link with any orientation. In these cases, the quandle is trivial of order 1 or 2, respectively.
\end{proposition}
{\bf Proof:} Suppose the fundamental quandle $Q(L)$ of the link $L$ is finite. Hence, every quotient of $Q(L)$ is finite and so $Q_n(L)$ is finite for all $n>1$. From the classification of links with finite $n$-quandles given in \cite{HS2},  the only non-trivial link for which  this is true is the Hopf link. A simple computation completes the proof.\hfill $\square$

Given a presentation of an $n$-quandle, one can try to systematically enumerate its elements and simultaneously produce a Cayley graph of the quandle. Such a method was described in a graph-theoretic fashion by Winker in \cite{WI}. The method is similar to the well-known Todd-Coxeter process for enumerating cosets of a subgroup of a group \cite{TC} and has been  extended to racks in \cite{HS3}. (A rack is more general than a quandle, requiring only axioms A2 and A3.) We provide a brief description of Winker's method applied to the $n$-quandle of a link since it will be used extensively in this paper.  Suppose $Q_n(L)$ is presented as
$$Q_n(L)=\langle x_1, x_2, \dots, x_g \, |\, x_{j_1}^{w_1}=x_{k_1}, \dots, x_{j_r}^{w_r}=x_{k_r} \rangle_n,$$
where each $w_i$ is a word in the free group on $\{x_1,\dots , x_g\}$. Throughout this paper presentations of $n$-quandles will not explicitly list the $n$-quandle relations $x_i^{x_j^n}=x_i$ (nor the relations given by the quandle axioms) although we are implicitly assuming they hold. To avoid confusion we append the subscript $n$ to presentations of $n$-quandles.

If $y$ is any element of the quandle, then it follows from the relation $x_{j_i}^{w_i}=x_{k_i}$ and Lemma~\ref{leftassoc} that $y^{\overline{w}_i x_{j_i}w_i}=y^{x_{k_i}}$, and so
$$y^{\overline{w}_i x_{j_i} w_i \overline{x}_{k_i}}=y.$$

Winker calls this relation the {\it secondary relation}  associated to the {\it primary relation} $x_{j_i}^{w_i}=x_{k_i}$. He also considers relations of the form $y^{ x_j^n }=y$ for all $y$ and $1 \le j \le g$. These relations are equivalent to the secondary relations of the $n$-quandle relations. In order to see this, notice that the secondary relation of the $n$-quandle relation $x_k^{x_j^n}=x_k$ is $y=y^{\bar x_j^n x_k x_j^n \bar x_k}$ for all elements $y$.  Now given any $z$, if we let $y=z^{x_j^n}$ in this secondary relation, then $z^{x_j^n}=z^{x_j^n \bar x_j^n x_k x_j^n \bar x_k} =  z^{x_k x_j^n \bar x_k}$. Hence, $z^{x_j^n x_k}=z^{x_k x_j^n}$ for all elements $z$. In a similar manner we find $z^{x_j^n \bar x_k}=z^{\bar x_k x_j^n}$ for all $z$.  Now given any $y$ we have $y=x_i^w$ for some $1 \le i \le g$ and, by what we just observed, we can commute $x_j^n$ with $w$ in the exponent of $x_i$. Therefore,
$$y^{x_j^n} = x_i^{w x_j^n} = x_i^{x_j^n w} = x_i^w =y.$$
Conversely, if $y^{x_j^n} =y$ for all $y$, then clearly $x_i^{x_j^n} =x_i$ as well.  

Winker's method now proceeds to build the Cayley graph associated to the presentation as follows:
  
\begin{enumerate}
\item Begin with $g$ vertices labeled $x_1,x_2, \dots, x_g$ and numbered $1,2, \dots,g$. 
\item Add an oriented loop at each vertex $x_i$ and label it $x_i$. (This encodes the axiom A1.)
\item For each value of $i$ from $1$ to $r$, {\em trace} the primary relation $x_{j_i}^{w_i}=x_{k_i}$ by introducing new vertices and oriented edges as necessary to create an oriented path from $x_{j_i}$ to $x_{k_i}$ given by $w_i$. Consecutively number (starting from $g+1$) new vertices in the order they are introduced.  Edges are labelled with their corresponding generator and oriented to indicate whether $x_i$ or $\overline x_i$ was traversed. 
\item Tracing a relation may introduce edges with the same label and same orientation into or out of a shared vertex. We identify all such edges, possibly leading to other identifications. This process is called {\it collapsing} and all collapsing is carried out before tracing the next relation. 
\item Proceeding in order through the vertices, trace and collapse each $n$-quandle relation $y^{x_j^n}=y$  and each secondary relation (in order). All of these relations are traced and collapsed at a vertex before proceeding to the next vertex.
\end{enumerate}

The method will terminate in a finite graph if and only if the $n$-quandle is finite. The reader is referred to Winker~\cite{WI} and \cite{HS3}  for more details.

Associated to every quandle $Q$ is its automorphism group $\text{Aut} (Q)$. The {\it inner} automorphism group of $Q$, denoted by $\text{Inn}(Q)$,  is the  normal subgroup of $\text{Aut}(Q)$ generated by the point symmetries $S_x$.  The {\it transvection} group of $Q$, denoted by $\text{Trans}(Q)$, is the subgroup of $\text{Inn}(Q)$ generated by all products $S_x S_y^{-1}$. The subgroup $\text{Trans}(Q)$ is normal in both $\text{Aut} (Q)$ and $\text{Inn}(Q)$. Moreover, $\text{Trans}(Q)$ is abelian if and only if 
\begin{equation}\label{abelian}(x^y)^{(z^w)}=(x^z)^{(y^w)}\end{equation}
for all elements $x, y, z, w\in Q$. Quandles that satisfy the property \eqref{abelian} are called {\it medial} or {\it abelian}. See \cite{JO} for more details. 

Some results in this paper were obtained using the RIG package for GAP.  The Cayley graph of a finite quandle $Q$ can be used to produce the operation table for $\rhd$ which is encoded in a matrix $M_Q$. In RIG, a rack (or quandle) can then be defined using the command {\sc Rack}$(M_Q)$. Once the quandle is entered into RIG, the built-in commands {\sc AutomorphismGroup}, {\sc InnerGroup}, and {\sc TransvectionsGroup} will compute $\text{Aut} (Q)$, $\text{Inn} (Q)$, and $\text{Trans} (Q)$, respectively.  Finally, the GAP command {\sc StructureDescription} will determine the structure of the group, such as $\mathbb Z_2 \times S_4$.  No additional special code is required to reproduce our results.

\section{Two-Bridge Links}\label{2bridge}

\begin{figure}[htbp]
\centerline{\includegraphics[height=1.5in]{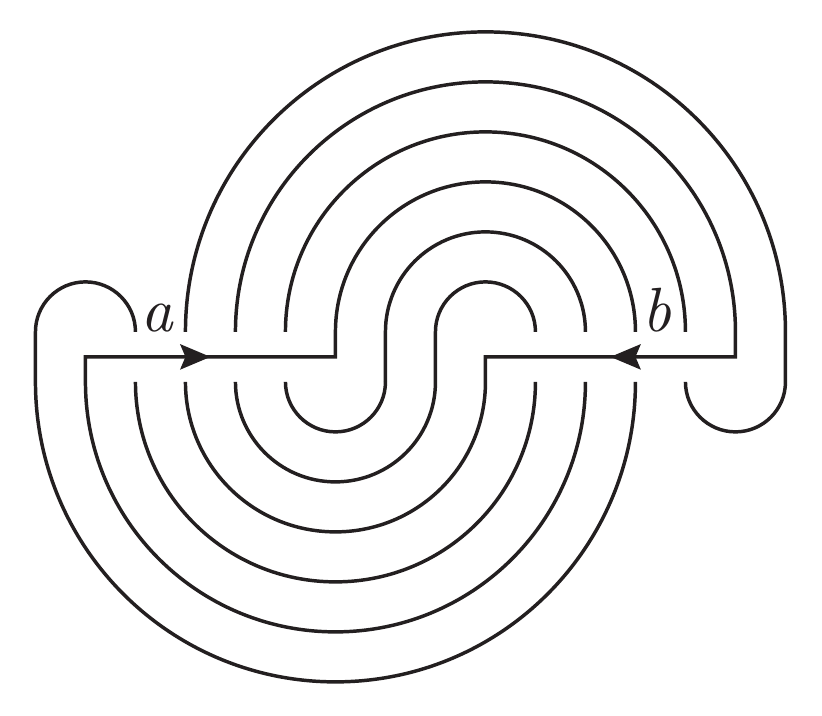}}
\caption{The Schubert normal form of the two-bridge link $L_{3/5}$.}
\label{L35}
\end{figure}

In this section we consider the involutory quandle of the non-trivial 2-bridge  link $L_{p/q}$. Because the 2-fold cyclic cover of $S^3$ branched over a non-trivial two-bridge link is a lens space with finite fundamental group, these quandles are finite.  We can easily find a presentation for the fundamental quandle of $L_{p/q}$, where $\gcd(p, q) = 1$ and $0 < p < q$ using the Schubert normal form (see, for example, \cite{BZ, KA}). As an example, consider the Schubert normal form of the figure-eight  knot $L_{3/5}$ with orientation shown in Figure~\ref{L35}.  The fundamental quandle for this knot has the presentation
$$Q(L_{3/5}) = \langle a, b \mid a^{b\bar{a}\bar{b}a} = b,\ b^{a\bar{b}\bar{a}b} = a \rangle,$$
where the generators $a$ and $b$ are the arcs of the two bridges.
 
For the involutory quandle, we may ignore the orientation of the link.  Thus, the involutory quandle of $L_{3/5}$ (with any orientation) is  presented by
$$Q_2(L_{3/5}) =  \langle a, b \mid a^{baba} = b,\ b^{abab} = a \rangle_2 = \langle a, b \mid a^{baba} = b \rangle_2,$$
where the two presentations are equivalent because, in the involutory quandle, the second relation is equivalent to the first relation.  In general, the presentation of $Q_2(L_{p/q})$ depends only on the number of undercrossings along one bridge in the Schubert normal form, which is equal to $q-1$. We have two cases, depending on whether $q$ is odd or even (i.e., whether $L_{p/q}$ is a knot or a link, respectively).

\begin{figure}[htbp]
\centerline{\includegraphics[width=3.0in]{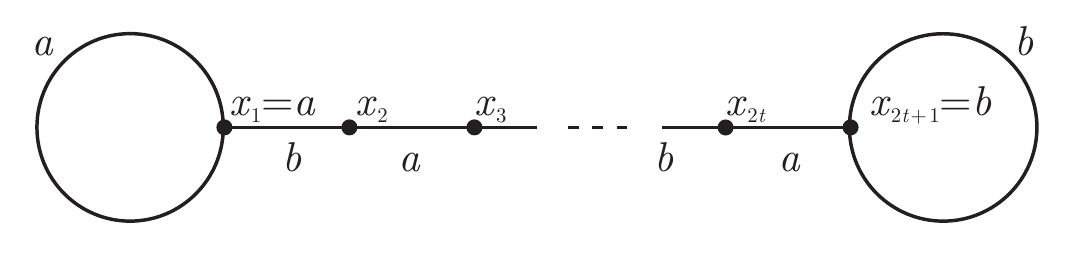}}
\caption{The Cayley graph for $Q_2(L_{p/q})$ with $q = 2t+1$.}
\label{Cayley2odd}
\end{figure}

If $q = 2t+1$ is odd, then a presentation of the involutory quandle is
$$Q_2(L_{p/q}) = \langle a, b \mid a^{(ba)^t} = b \rangle_2.$$
The secondary relation associated to the primary relation $a^{(ba)^t} = b$,  is $x^{(ab)^{2t+1}} = x$.  Tracing the primary relation by Winker's method gives the diagram in Figure~\ref{Cayley2odd} (where the vertex $x_1$ represents the element $a$ and $x_{2t+1}$ the element $b$). 

Notice from Figure~\ref{Cayley2odd} that 
\begin{equation}\label{abpath} x_i^{(ab)^t} =\left\{ \begin{array}{ll} x_1, & i=2t+1 \\ x_{2t+2-i+(-1)^i}, & 
\mbox{$1\le i < 2t+1$.}
\end{array} \right. 
\end{equation}

It is now straightforward to check that the secondary relation is satisfied at each vertex.  For example, if $1 \le i < 2t+1$ is odd, then by (\ref{abpath}) we have
$$x_i^{(ab)^{2t+1}} = x_i^{(ab)^t(ab)^tab}=x_{2t+1-i}^{(ab)^tab}=x_{i+2}^{ab}=x_i.$$
The remaining cases are similar. Therefore, Figure~\ref{Cayley2odd} gives the Cayley graph of $Q_2(L_{p/q})$.  The number of elements in the quandle is $2t+1 = q$.

\begin{figure}[htbp]
\centerline{\includegraphics[width=3.0in]{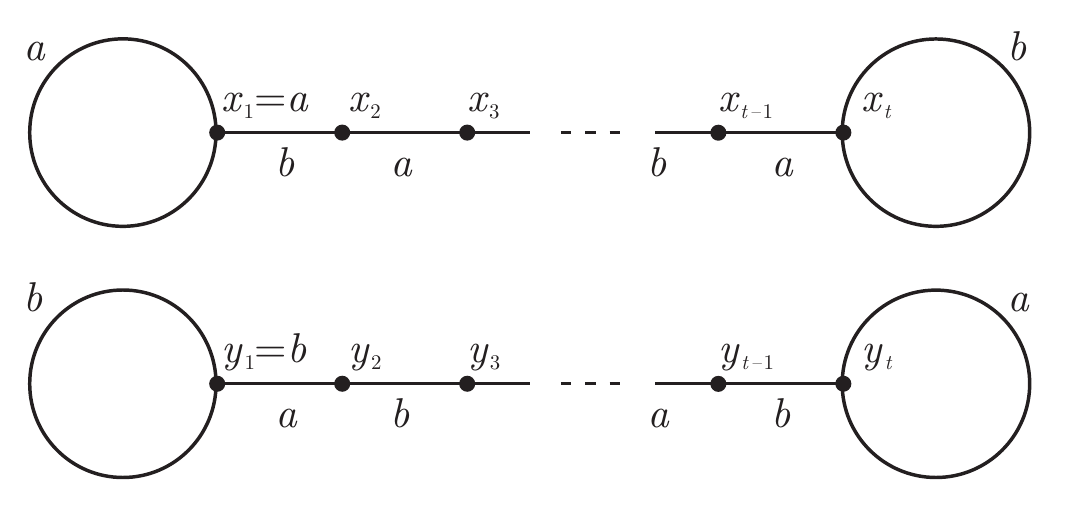}}
\caption{The Cayley graph for $Q_2(L_{p/q})$ with $q = 2t$ and $t$ odd.}
\label{Cayley2even}
\end{figure}

If $q = 2t$ is even, then a presentation of the involutory quandle is
$$Q_2(L_{p/q}) = \langle a, b \mid a^{(ba)^{t-1}b} = a,\ b^{(ab)^{t-1}a} = b \rangle_2.$$
In this case there are two primary relations.  However, the secondary relations associated to the two primary relations are the same: $x^{(ab)^{2t}} = x$.  Tracing the primary relations gives the diagram in Figure~\ref{Cayley2even} when $t$ is odd (where the vertex $x_1$ is the element $a$ and $y_1$ is the element $b$).  It is once again straightforward to check that the secondary relation is satisfied at each vertex, so this diagram is the Cayley graph of the involutory quandle.  Notice that there are two {\em algebraic} components of the quandle, each given by a connected component of the Cayley graph. This reflects the fact that $L_{p/q}$ with $q$ even is a 2-component link. In general, the number of link components is equal to the number of algebraic components of its quandle (see \cite{JO}).  The number of elements in each algebraic component of $Q_2(L_{p/q})$ is $t$, and so the total number of elements is $2t = q$. A similar analysis applies to $t$ even.

In summary, for any non-trivial two-bridge knot or link, the involutory quandle $Q_2(L_{p/q})$ has order $q$. In fact,  $Q_2(L_{p/q})$ is isomorphic to the  {\em dihedral quandle} which is the set $R_q=\mathbb Z/q \mathbb Z$ with quandle operations defined by $i^j=i^{\bar j} = 2j-i\  (\mbox{mod }q)$. 

\begin{proposition} If $L_{p/q}$ is a non-trivial two-bridge knot or link, then $Q_2(L_{p/q})$ is isomorphic to the dihedral quandle $R_q$.
\end{proposition} 

{\bf Proof.} There are two cases depending on whether $q$ is even or odd. We present only the odd case; the remaining case is similar.  Assume $q = 2 t+1$. From above we have the presentation
$$Q_2(L_{p/q}) = \langle a, b \mid a^{(ba)^t} = b \rangle_2.$$
Define a map $\phi: Q_2(L_{p/q}) \rightarrow R_q$ by defining $\phi$ on the generators by $\phi(a)=0$ and $\phi(b)=1$ and then extending the map over $Q_2(L_{p/q})$ using the quandle operation.  (That is, if $\phi(x)=i$ and $\phi(y)=j$, then define $\phi(x^y) = i^j$.) To verify that $\phi$ is well-defined, it suffices to check that the image of the relation $a^{(ba)^t} = b$ is satisfied in $R_q$.  Notice that if $\phi(x)=i$, then 
$$\phi(x^{ba})=  \left(i^1 \right)^0 = i-2\ (\mbox{mod }q).$$
It then follows that 
$$\phi \left( a^{(ba)^t} \right) = -2t = 1 = \phi(b) \ (\mbox{mod }q).$$
The map $\phi$ is a homomorphism by definition. Moreover, it is surjective because $\phi(a)=0$ and $\phi(b)=1$ generate $R_q$. This follows inductively from the observation that $i^{i+1} = i+2\ (\mbox{mod }q)$ for all $i \in R_q$. Finally, since $|Q_2(L_{p/q})|=|R_q|=q$, we have that $\phi$ is an isomorphism. \hfill $\Box$

In \cite{EM}, Elhamdadi, MacQuarrie, and Restrepo prove that $\mbox{Aut}(R_n)$ is the semi-direct product $\mathbb Z_n \rtimes \mathbb Z_n^*$ and $\mbox{Inn}(R_n)$ is the dihedral group $D_n$ of order $2n$ if $n$ is odd and $D_{n/2}$ if $n$ is even. Here $\mathbb{Z}_n^*$ is the multiplicative group of units in $\mathbb{Z}_n$, and the semidirect product is given by the homomorphism $\phi:\mathbb Z_n \rightarrow {\rm Aut}(\mathbb Z_n)$ defined by $\phi(a)(x) = ax$. Thus, we obtain the results given in Table~\ref{automorphism groups of 2-bridge links}. We further note that  $\textrm{Trans}(Q_2(L_{p/q}))$ is generated by a rotation in ${\rm Inn}(Q_2(L_{p/q})).$ 

\begin{table}[h]
\begin{center}
\renewcommand{\arraystretch}{1.5}
\begin{tabular}{|c||c|c|c|c|}\hline
$Q_n$ & $|Q_n|$ & \textrm{Aut} & \textrm{Inn} & \textrm{Trans} \\ \hline
$Q_2(L_{p/q})$ & $q$ & $\mathbb Z_q \rtimes \mathbb Z_q^*$ & $D_{q/ \gcd(2,q)} $ & $\mathbb Z_{q/ \gcd(2,q)}$ \\ \hline
\end{tabular}
\end{center}
\caption{The automorphism groups of $Q_2(L_{p/q})$. }
\label{automorphism groups of 2-bridge links}
\end{table}

\section{Torus Links}

As mentioned in the introduction, the $n$-quandle of a link $L$ is finite if and only if the $n$-fold cyclic branched cover of $S^3$ branched over $L$ has finite fundamental group. These  covers, in the case where $L$ is the torus link $T_{p,q}$, were classified by Milnor in \cite{M} and as a result, $Q_n(T_{p,q})$ is finite if and only if $\frac{1}{p} + \frac{1}{q} + \frac{1}{n} >1$. The inequality holds for the following values of $p$, $q$, and $n$ given in the following table. 

\begin{table}[h]
\begin{center}
\renewcommand{\arraystretch}{1.5}
\begin{tabular}{|c||c|c|c|c|c|c|c|c|}\hline
$(p,q)$ &$(2,2)$ &$(2,3)$ &$(2,4)$ &$(2,5)$ &$(2,q),q>5$&$(3,3)$&$(3,4)$ &$(3,5)$ \\ \hline
$n$& $n>1$& 2, 3, 4, 5& 2, 3& 2, 3& 2& 2& 2& 2\\ \hline
\end{tabular}
\end{center}
\caption{Values of $p$, $q$, and $n$ for which $Q_n(T_{p,q})$ is finite.}
\label{pqnforfiniteQ}
\end{table}

\begin{figure}[htbp]
$$\begin{array}{ccc}
\includegraphics[height=1.1in,trim = 0 -30pt 0 0, clip]{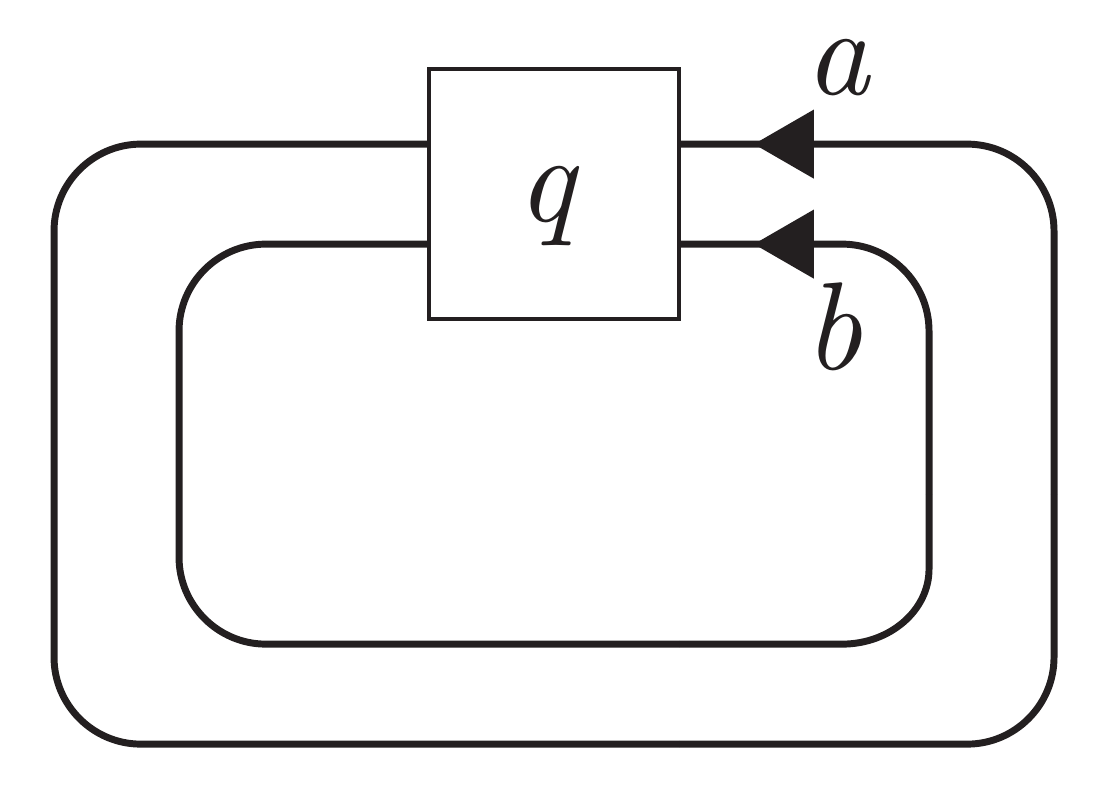} & & \includegraphics[height=1.25in,trim = 0 25pt 0 0, clip]{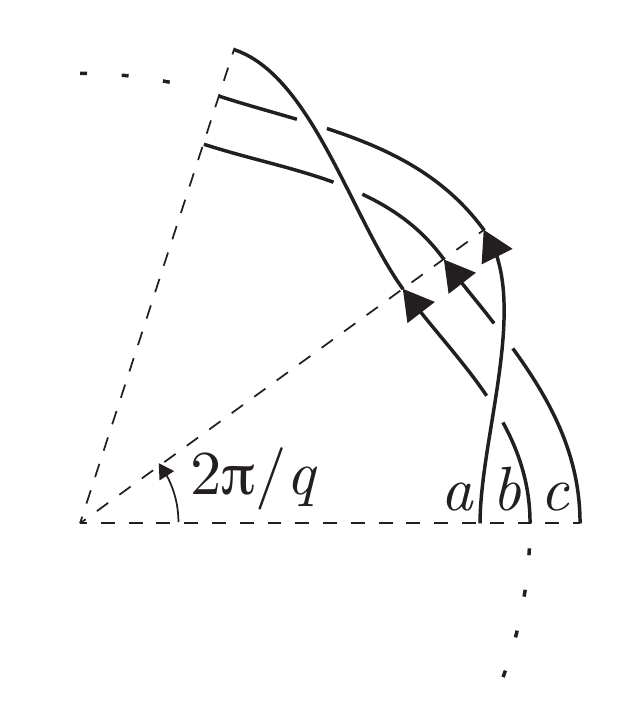} \\
T_{2,q} &\hspace{.25in} & T_{3,q}
\end{array}$$
\caption{The oriented links $T_{2,q}$ and $T_{3,q}$ with choice of generators.}
\label{torus diagrams}
\end{figure}

Since $T_{2,q}$ is a two-bridge link, the links $T_{2,q}$ with $n=2$ in Table~\ref{pqnforfiniteQ} were considered in the previous section. For each of the remaining cases, it is a simple matter to derive a presentation from a link diagram and then employ Winker's method to create a Cayley graph of the associated quandle. Diagrams for $T_{2,q}$ and $T_{3, q}$ with choice of orientation and quandle generators are shown in Figure~\ref{torus diagrams}. In the case of $T_{2,q}$, the box labeled $q$ contains $q$ right handed half twists.  The remainder of this section lists the results. 
Note that only in the case of the two-component link $T_{2,4}$ with $n=3$ does orientation matter. We denote by $T_{2,4}^{+-}$ the oriented link obtained from $T_{2,4}$ by reversing the orientation of the second component (the one labeled $b$).

\begin{figure}[htbp]
$$\begin{array}{cccc}
\includegraphics[height=1.0in]{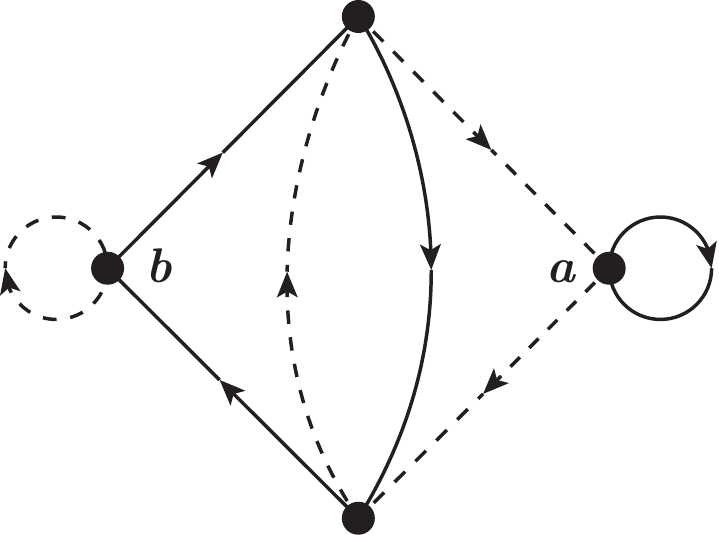} & \includegraphics[height=1.0in]{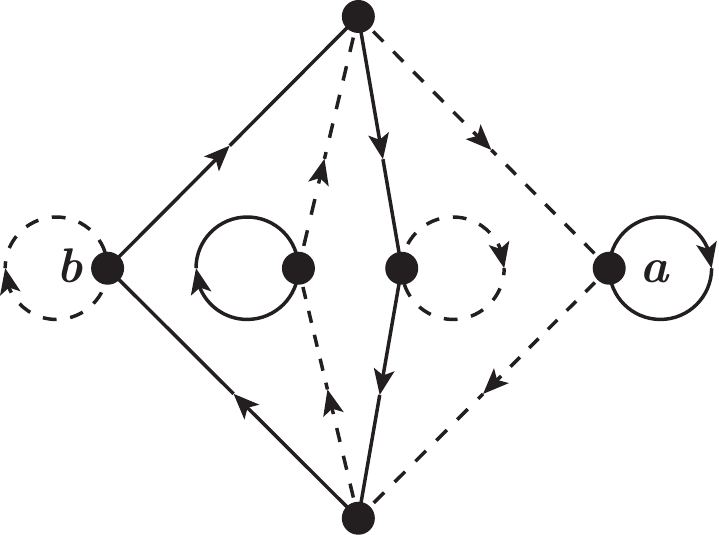} &
\multicolumn{2}{c}{\includegraphics[height=1.0in]{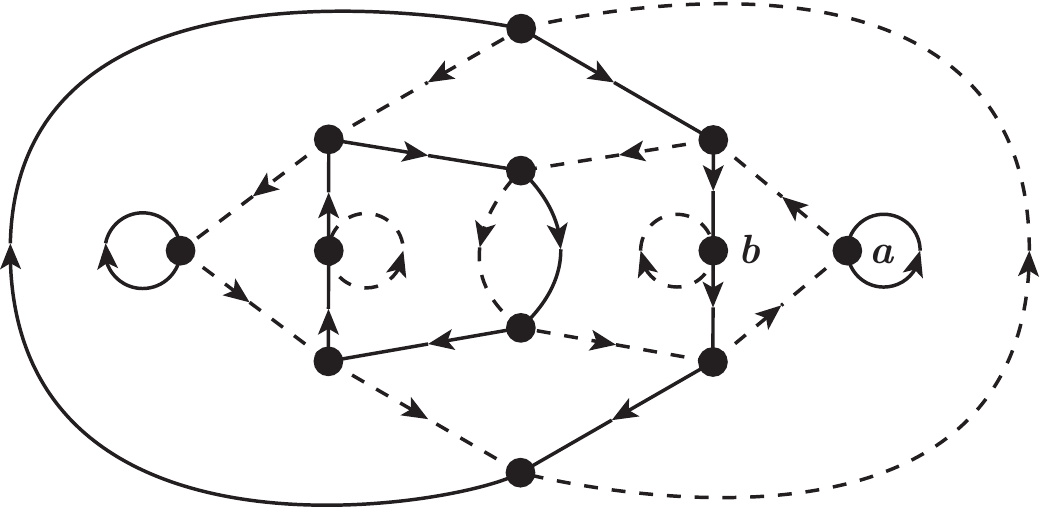}} \\
\\
Q_3(T_{2,3}) & Q_4(T_{2,3}) & \multicolumn{2}{c}{Q_5(T_{2,3})} \\
\\
\multicolumn{2}{c}{\includegraphics[height=1.0in]{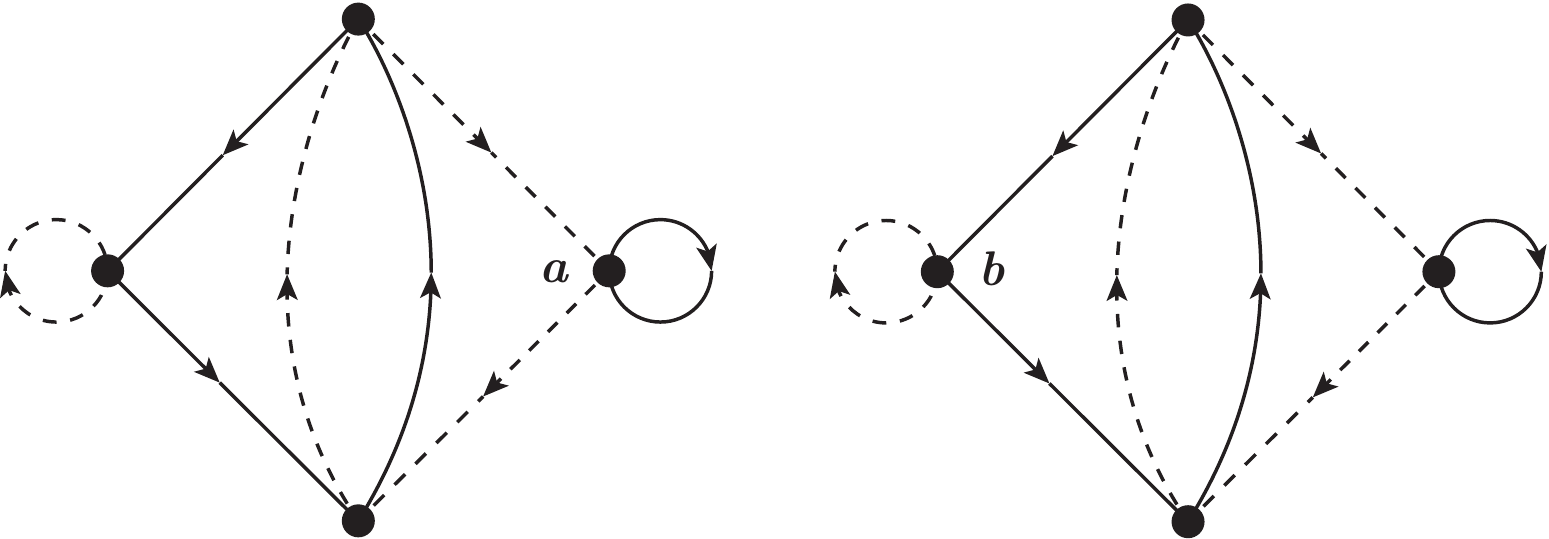}} & \multicolumn{2}{c}{\includegraphics[height=1.0in]{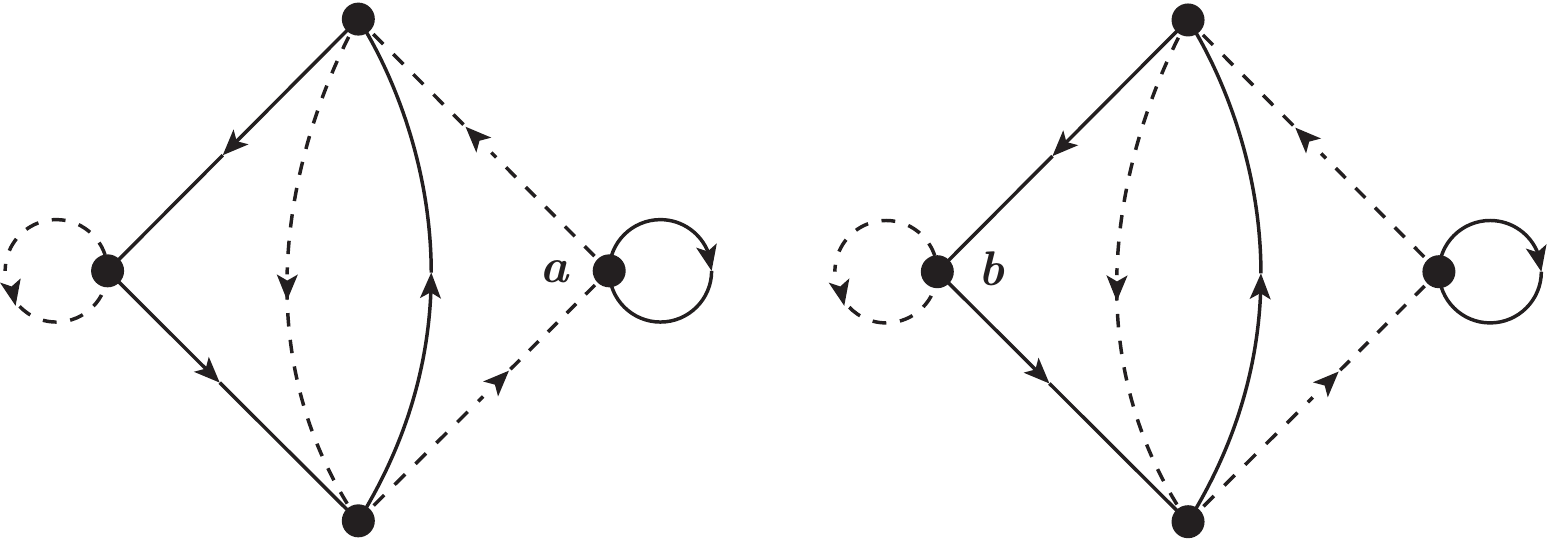}} 
\\
\\
\multicolumn{2}{c}{Q_3(T_{2,4})} & \multicolumn{2}{c}{Q_3(T_{2,4}^{+-})}
\\
\\
\multicolumn{4}{c}{\includegraphics[height=1.5in]{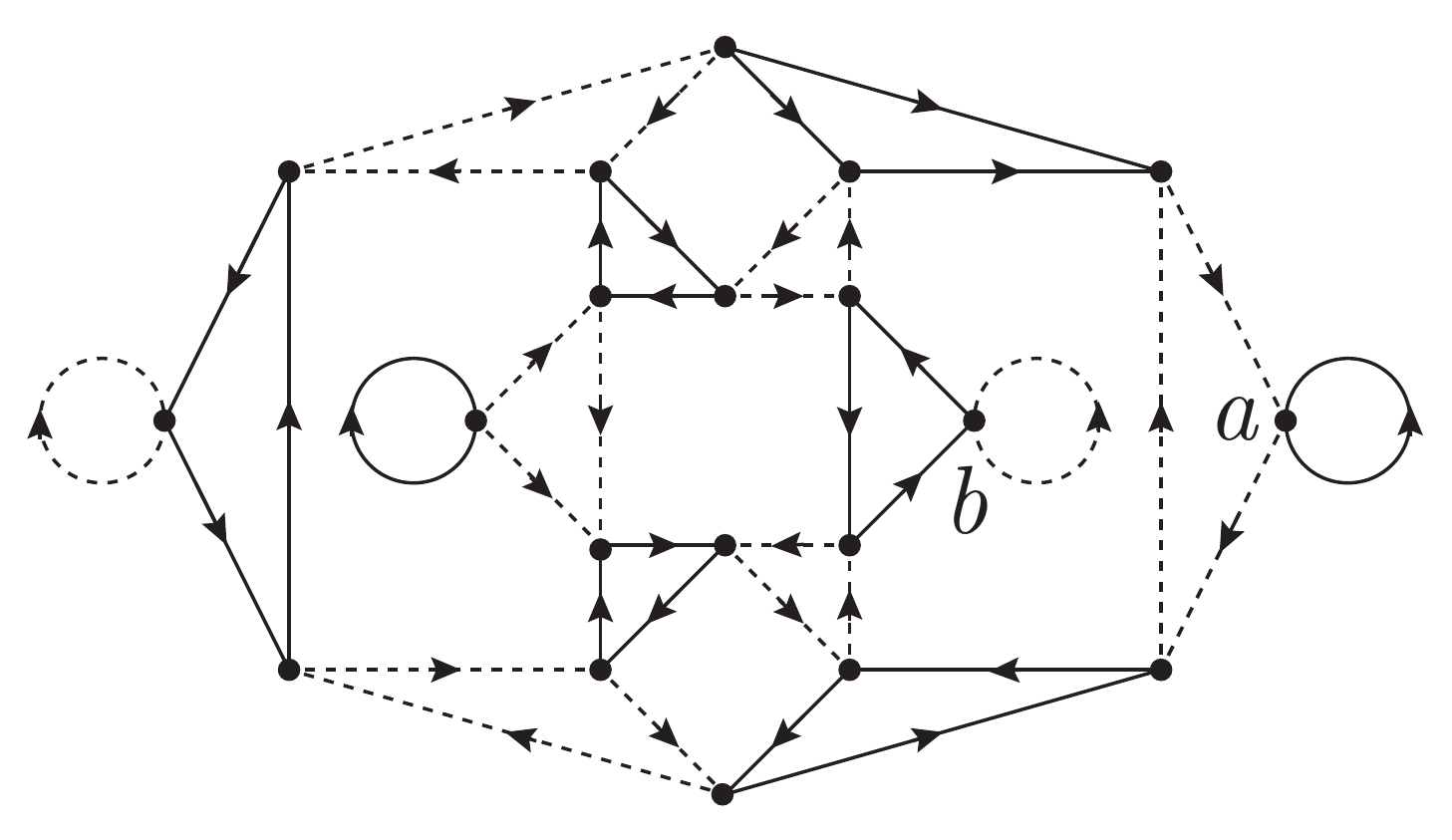}}
\\
\multicolumn{4}{c}{Q_3(T_{2,5})}
\end{array}$$
\caption{Cayley graphs of $Q_n(T_{2,q})$ with $n > 2$.}
\label{torus Cayley graphs part 1}
\end{figure}

\begin{figure}[htbp]
$$\begin{array}{ccccc}
\includegraphics[height=1.25in,trim= 0pt -30pt 0pt 0pt,clip]{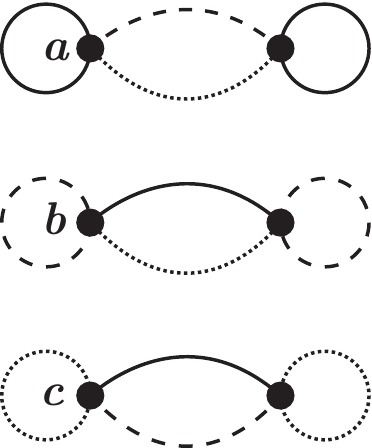} & \hspace{.2in} &   
\includegraphics[height=1.5in]{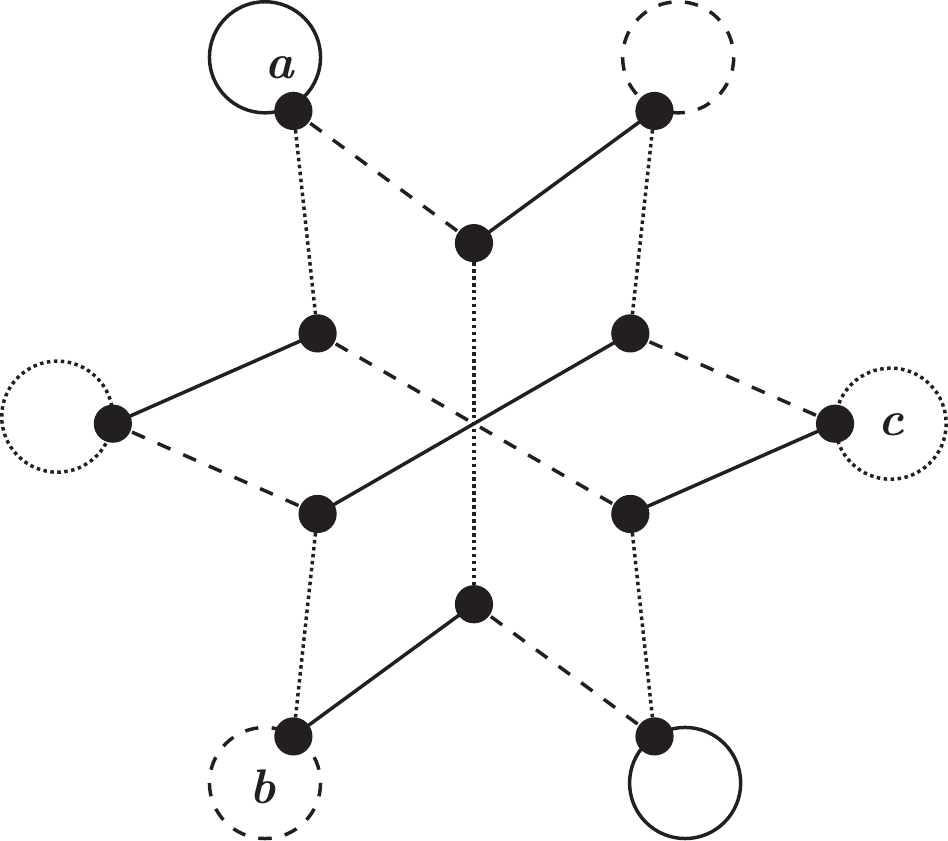} & \hspace{.2in} & 
\includegraphics[height=1.5in]{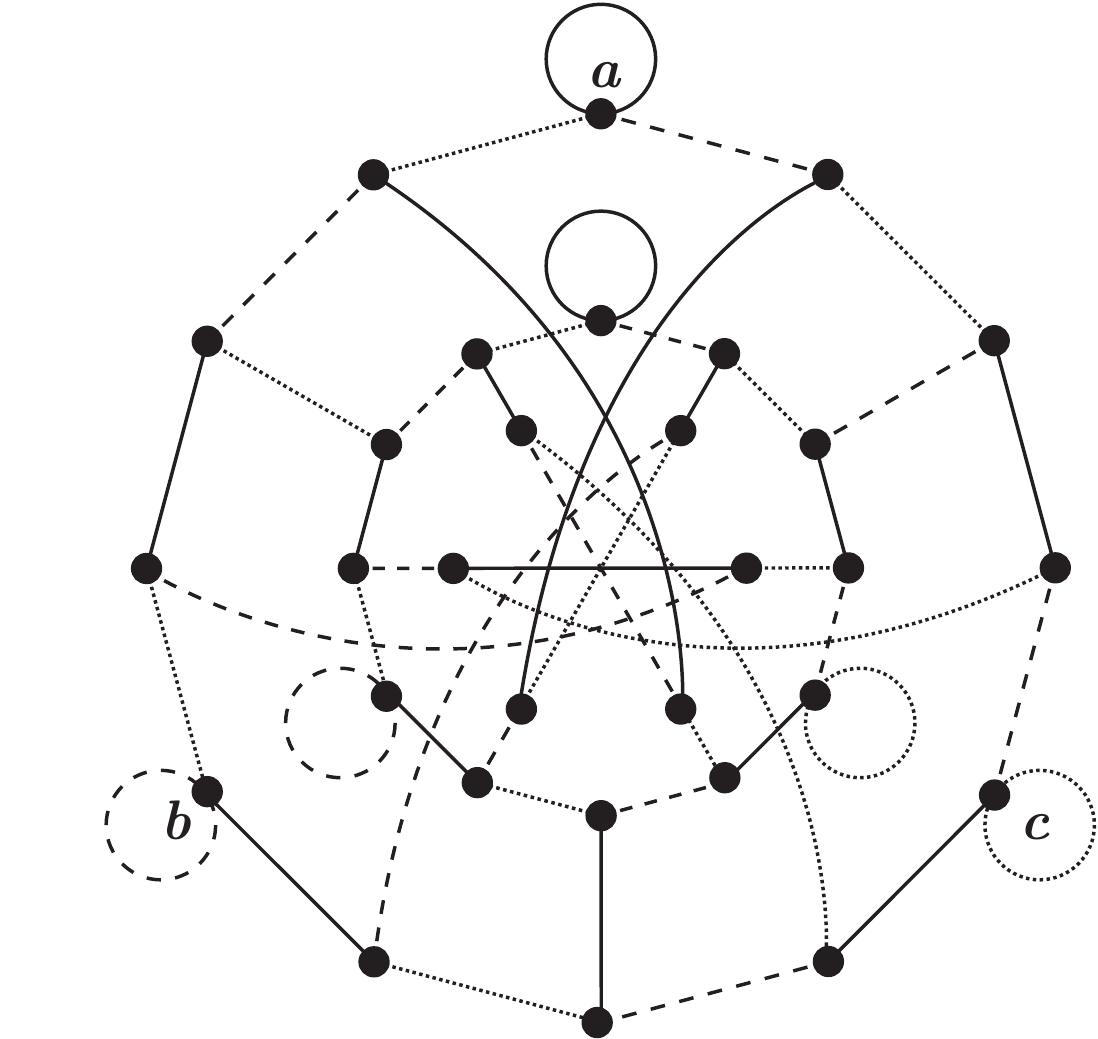} \\
\\
Q_2(T_{3,3}) & & Q_2(T_{3,4}) & & Q_2(T_{3,5})
\end{array}$$
\caption{Cayley graphs of $Q_2(T_{3,q})$.}
\label{torus Cayley graphs part 2}
\end{figure}

Cayley graphs of these quandles are given in Figures~\ref{torus Cayley graphs part 1} and \ref{torus Cayley graphs part 2}. In these graphs, solid, dashed, and dotted edges correspond to the generators $a$, $b$, and $c$, respectively. The automorphism groups were computed using GAP and are given in Table~\ref{torus quandle automorphism groups}.

\begin{table}
\begin{center}
\renewcommand{\arraystretch}{1.5}
\begin{tabular}{|c||c|c|c|c|}
\hline
$Q_n$ & $|Q_n|$ & \textrm{Aut} & \textrm{Inn} & \textrm{Trans} \\ \hline
$Q_2(T_{2,q})$ & $q$ & $\mathbb Z_q \rtimes \mathbb Z_q^*$ & $D_{q/\gcd(2,q)}$ & $\mathbb Z_{q/\gcd(2,q)}$ \\ \hline
$Q_3(T_{2,3})$ & 4 & $A_4$ & $A_4$ & $\Z_2 \times \Z_2$ \\ \hline
$Q_4(T_{2,3})$ & 6 & $S_4$ & $S_4$ & $A_4$ \\ \hline
$Q_5(T_{2,3})$ & 12 & $A_5$ & $A_5$ & $A_5$ \\ \hline
$Q_3(T_{2,4})$ & 8 & $S_4$ & $A_4$ & $A_4$ \\ \hline
$Q_3(T_{2,4}^{+-})$ & 8 & $\Z_2 \times A_4$ & $A_4$ & $\Z_2 \times \Z_2$ \\ \hline
$Q_3(T_{2,5})$ & 20 & $S_5$ & $A_5$ & $A_5$ \\ \hline
$Q_2(T_{3,3})$ & 6 & $\Z_2 \times S_4$ & $\Z_2 \times \Z_2$ & $\Z_2 \times \Z_2$  \\ \hline
$Q_2(T_{3,4})$ & 12 & $\Z_2 \times S_4$ & $A_4$ & $A_4$  \\ \hline
$Q_2(T_{3,5})$ & 30 & $\Z_2 \times S_5$ & $A_5$ & $A_5$  \\ \hline
\end{tabular}
\end{center}
\caption{Order and automorphism groups of finite torus link $n$-quandles.}
\label{torus quandle automorphism groups}
\end{table}

\section{Torus Links with Axis}

\begin{figure}[htbp]
$$\begin{array}{ccc}
\includegraphics[height=1.00in]{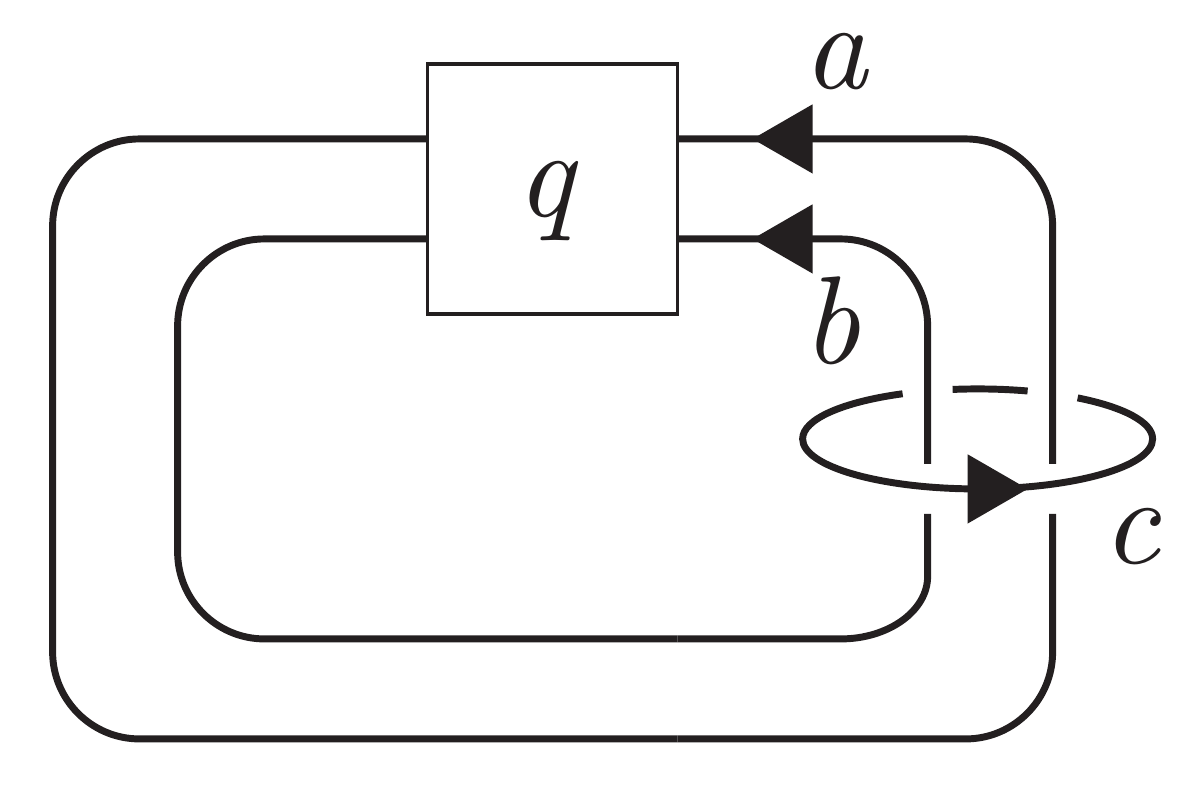} & &
\includegraphics[height=1.25in]{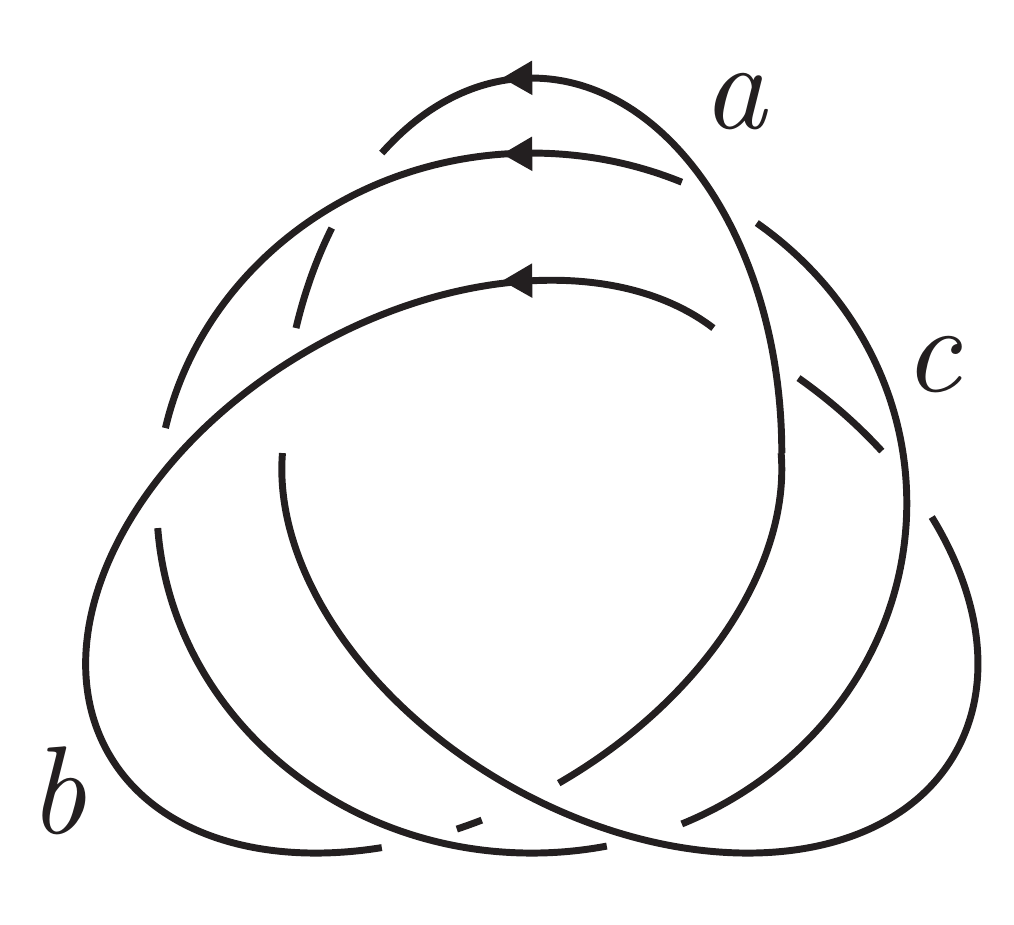}\\
T_{2,q}\cup A & & T_{2, 3}\cup B
\end{array}$$
\caption{The oriented links $T_{2,q}\cup A$ and $T_{2, 3}\cup B$ with choice of generators.}
\label{torus diagrams with axis}
\end{figure}

Given the torus link $T_{p,q}$ lying on the torus $F$ that separates $S^3$ into two solid tori, an {\em axis} of $T_{p,q}$ is the core of either solid torus.  In this section we consider the oriented torus links with axis shown in Figure~\ref{torus diagrams with axis}. As in the last section, $q$ represents $q$ right handed half twists. For each torus link $T_{2, q}$, we adjoin the axis $A$ with linking number $+2$ (as opposed to $+q$). In the case of the trefoil, we also include the oriented axis $B$ with linking number $+3$. Only the involutory quandles of these links are finite, hence the orientations are immaterial.  Since $T_{2,-q} \cup A = (T_{2,q} \cup A)^*$ and $T_{2,1}\cup A$ is the 2-bridge link $L_{1/4}$, we may further assume $q>1$. A choice of generators for each link is given in the figure.

\begin{figure}[htbp]
\centerline{\includegraphics*[height=1.5in]{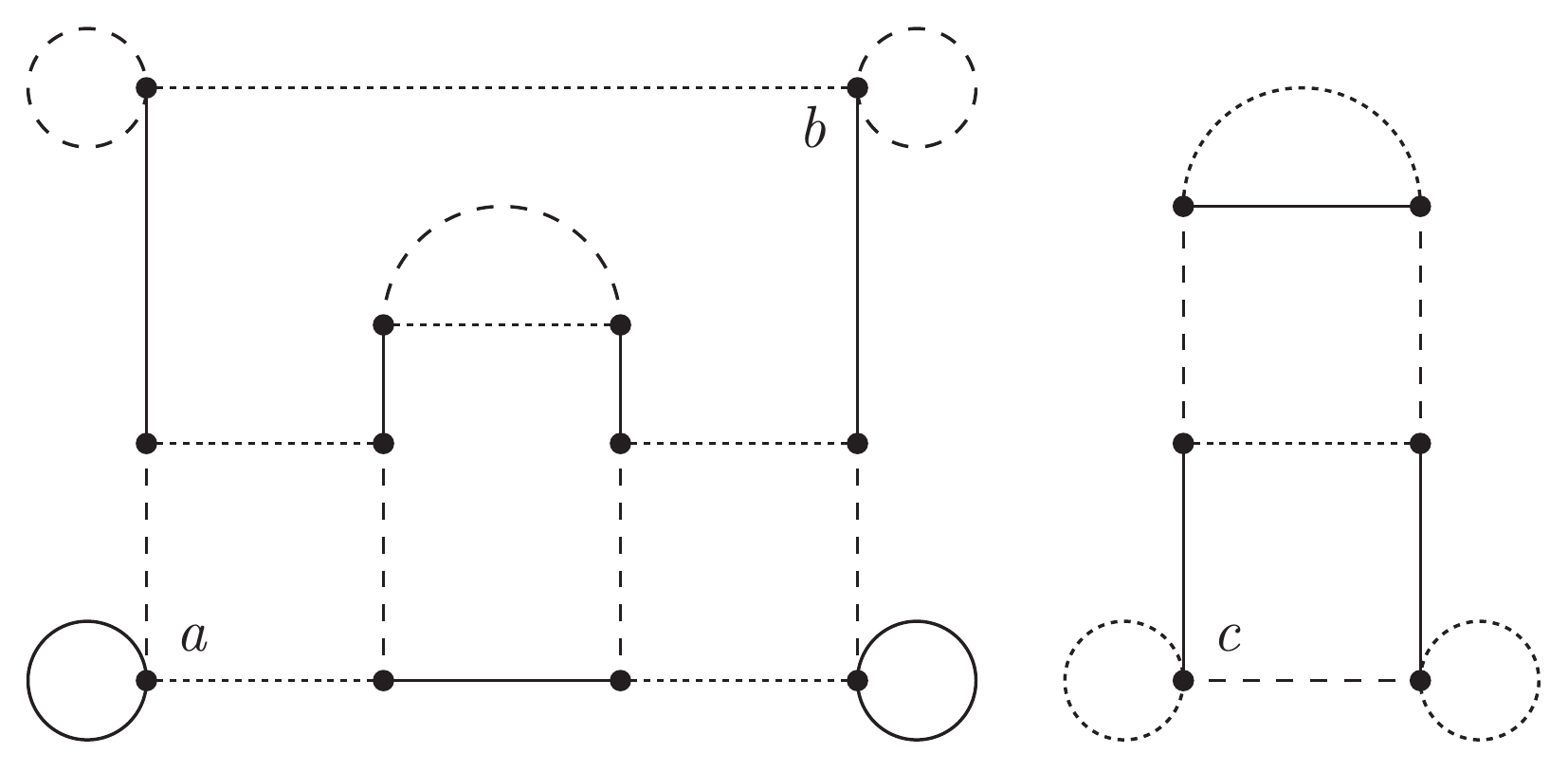}}
\caption{The Cayley graph of $Q_2(T_{2,3} \cup B)$. }
\label{T23UB cayley diagram}
\end{figure}

First, we consider the link $T_{2, 3}\cup B$. After deriving a presentation for $Q_2(T_{2,3} \cup B)$ from the diagram and employing Winker's method, we obtain the Cayley graph shown in Figure~\ref{T23UB cayley diagram}. Solid, dashed, and dotted edges correspond to the generators $a$, $b$, and $c$, respectively, in the figure. Using GAP to compute the automorphism groups of this quandle we find that ${\rm Aut}(Q_2(T_{2,3}\cup B)) \cong \Z_2 \times \Z_2 \times S_4$, ${\rm Inn}(Q_2(T_{2,3} \cup B)) \cong S_4$ and ${\rm Trans}(Q_2(T_{2,3} \cup B)) \cong S_4$.

Now we analyze the infinite family $T_{2,q}\cup A$. There are two cases depending on whether $q$ is even or odd. Assume $q = 2t+1$ with $t > 0$. A presentation of the involutory quandle in this case is
$$Q_2(T_{2,q} \cup A) = \langle a, b, c \mid c^{ab} = c,\ a^{(ba)^t bc} = b,\ b^{(ab)^t c} = a\rangle_2.$$
In order to make Winker's method simpler we change to the equivalent presentation
$$Q_2(T_{2,q} \cup A) = \langle a, b, c \mid \{c^{(ab)^i} = c \}_{i=1}^{2t},\ a^{c(ba)^t} = b,\ a^{(ba)^t c} = b,\ a^{cac}=a,\ b^{cbc}=b \rangle_2.$$
To see that the presentations are equivalent, we need to show that each set of relations can be derived from the other. Let $P_1$, $P_2$, and $P_3$ be the relations in the first presentation and $S$, $R_2$, $R_3$, $R_4$, and  $R_5$ be the relations in the second presentation in the order given. The relations $c^{(ab)^i} = c$ in the set $S$ follow from induction and relation $P_1$. The relation $R_2$ follows immediately from $P_3$.  The secondary relation associated to $P_1$ implies $x^{cba}=x^{bac}$ for all $x$ and thus $x^{c(ba)^j}=x^{(ba)^jc}$ for $j \ge 1$ follows by induction. Similarly, we have $x^{c(ab)^j}=x^{(ab)^jc}$ for $j \ge 1$.  The relation $R_2$ now follows from this observation and $P_3$.
The relation $R_4$ is then derived as follows:
$$a^{cac} \stackrel{P_3}{=} b^{(ab)^tac}= b^{b(ab)^t ac} = b^{(ba)^{t+1}c} = b^{c(ba)^{t+1}}\stackrel{P_2}{=} a^a = a.$$
The derivation of $R_5$ is similar:
$$b^{cbc} \stackrel{R_3}{=} a^{(ba)^tbc} \stackrel{P_2}{=} b.$$
In a similar manner, $P_1$, $P_2$, and $P_3$ can be derived from $S$, $R_2$, $R_3$, $R_4$, and  $R_5$.

\begin{figure}[htbp]
\centerline{\includegraphics*[height=1.25in]{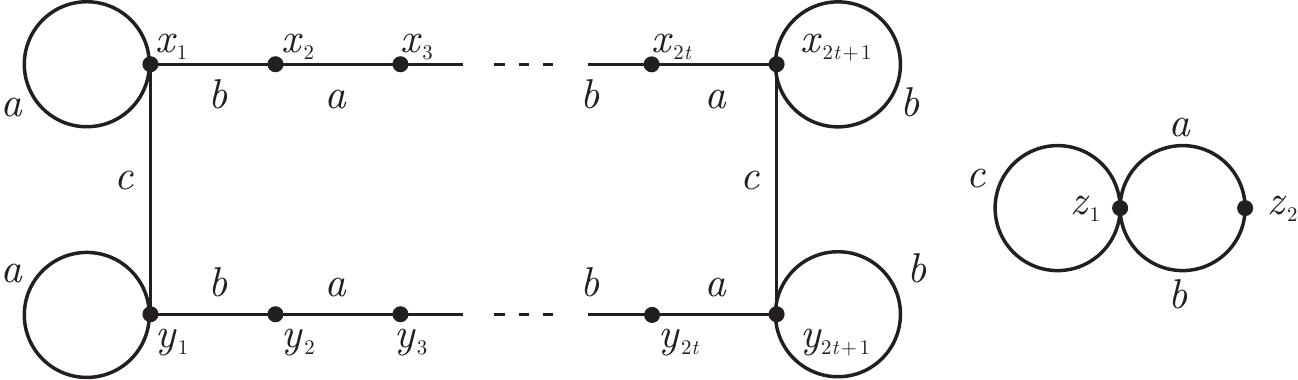}}
\caption{The primary relation graph for $Q_2(T_{2,q} \cup A)$ with $q=2t+1$.}
\label{Q2T2qAoddprm}
\end{figure}

We now proceed with the diagramming method. We begin with three vertices $x_1=a$, $y_{2t+1}=b$, and $z_1=c$ and loops labeled $a$, $b$, $c$, respectively, at each vertex. We next trace the $2t+4$ primary relations which introduces the vertices $z_2 ,y_1, \dots, y_{2t},x_2, \dots, x_{2t+1}$ in that order. Tracing the primary relations gives the graph shown in Figure~\ref{Q2T2qAoddprm}. 

\begin{figure}[htbp]
\vspace*{13pt}
\centerline{\includegraphics*[height=1.25in]{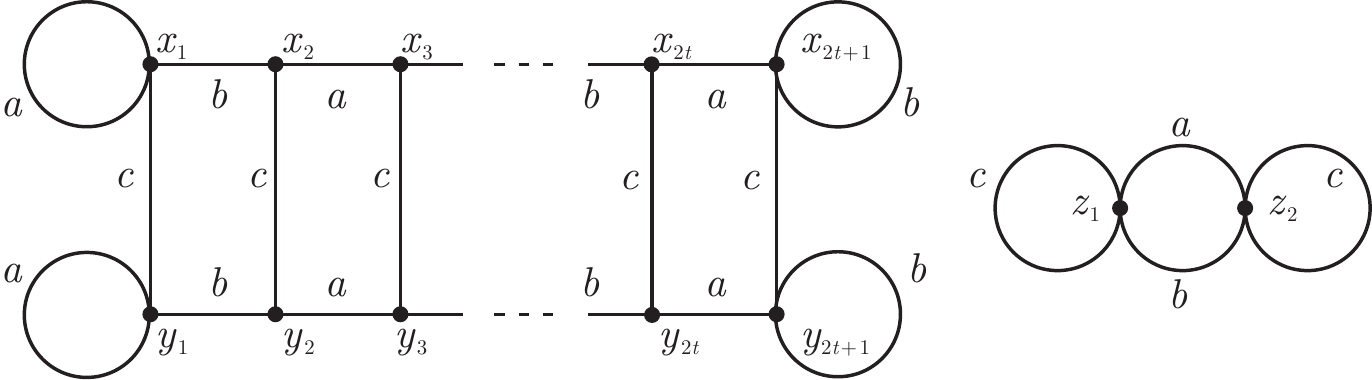}}
\caption{The Cayley graph for $Q_2(T_{2,q} \cup A)$ with $q=2t+1$.}
\label{Q2T2qJodd}
\end{figure}

Finally, we consider the secondary relations $x^{w_j}=x$ for $1 \le j \le 5$ where,
\begin{align*}
w_1 =& (ba)^ic(ab)^ic, \ 1 \le i \le 2t, \\
w_2 =& (ab)^t cac b(ab)^t, \\
w_3 =& c (ab)^{2t} a c b, \\
w_4 =& (ca)^4, \\
w_5 =& (cb)^4.
\end{align*}
We trace the secondary relations at each vertex in the order the vertices were introduced. Notice that the  formulas in (\ref{abpath}) from Section~\ref{2bridge} apply to both $x_i$ and $y_i$ in Figure~\ref{Q2T2qAoddprm}. Tracing the secondary relation $w_1$ for $1 \le i \le 2t$ at vertex $x_1$ introduces $2t-1$ edges labeled $c$ that connect the vertices $x_i$ to $y_i$. At this point we claim that all secondary relations are satisfied at vertices $x_i$ and $y_i$ for $1 \le i \le 2t+1$. There are several cases to consider. We will verify the relations $x_i^{w_2}=x_i$ and 
$x_i^{w_3}=x_i$ for $i \ne 1$ and odd and leave the remaining cases to the reader. Using (\ref{abpath}) and Figure~\ref{Q2T2qAoddprm} we have:
$$x_i^{(ab)^t cac b(ab)^t}= x_{2t+1-i}^{cac b(ab)^t} = x_{2t+3-i}^{(ab)^t} = x_i.$$
Similarly:
$$x_i^{c (ab)^{2t} a c b}= y_i^{(ab)^{2t} a c b} = y_{2t+1-i}^{(ab)^{t} a c b} = y_{i+2}^{a c b} = x_i.$$
Finally, consider tracing the secondary relations at vertex $z_1$. The relations $z_1^{w_1}=z_1$ for $1 \le i \le 2t$ are already satisfied. Tracing $z_1^{w_2}=z_1$ introduces a loop at $z_2$ labeled $c$. This gives the graph in Figure~\ref{Q2T2qJodd}. It is not hard to then verify that all remaining secondary relations are satisfied at $z_1$ and $z_2$. Therefore, Figure~\ref{Q2T2qJodd} is the Cayley graph of $Q_2(T_{2,q} \cup A)$ with $q$ odd. From this we see that $Q_2(T_{2,q} \cup A)$ with $q$ odd has order $2q+2$. Notice that the graph has two connected components which was expected since $T_{2,q} \cup A$ is a link of 2 components.

\begin{figure}[htbp]
\centerline{\includegraphics*[height=2.75in]{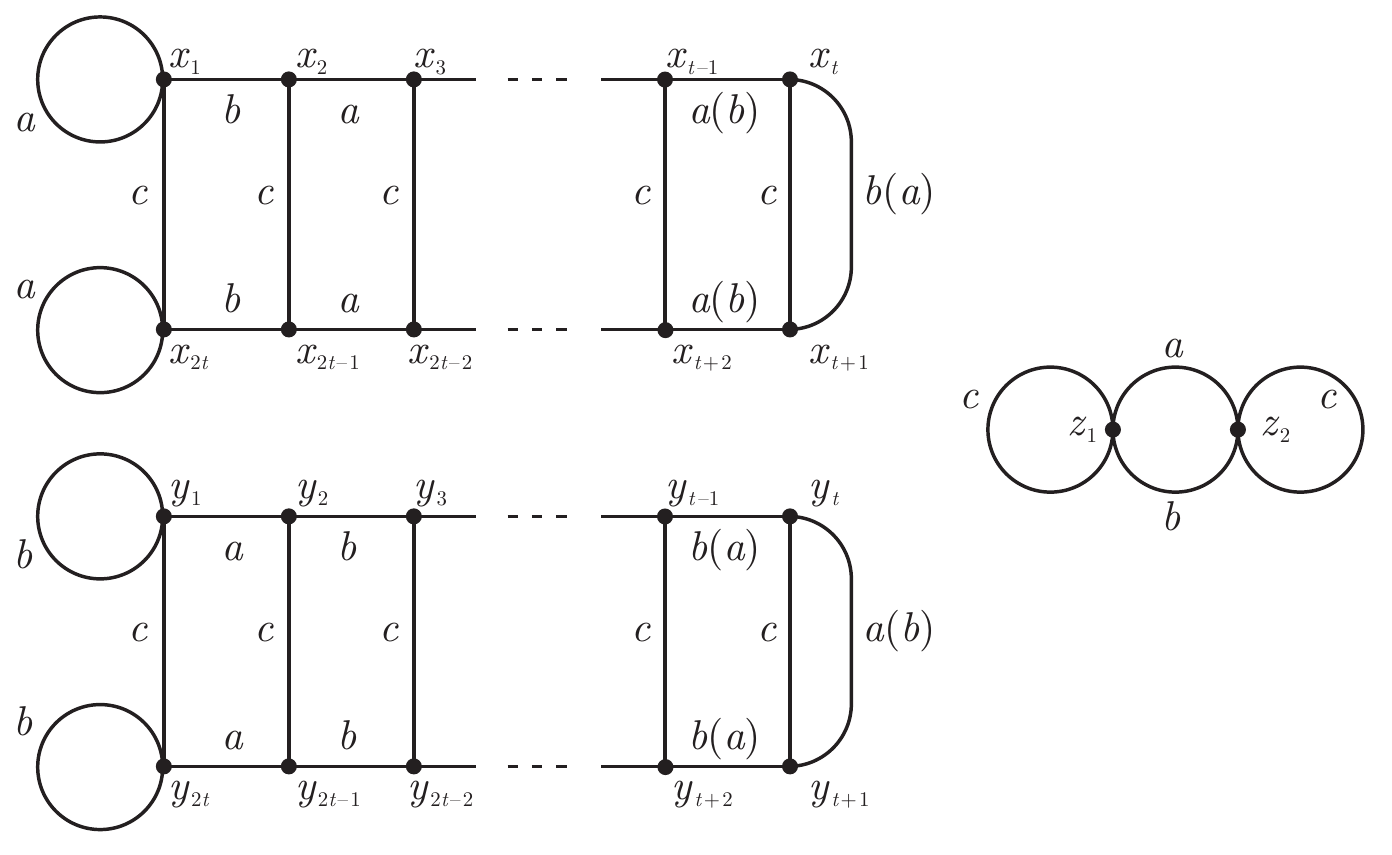}}
\caption{The Cayley graph for $Q_2(T_{2,q} \cup A)$ with $q=2t$.}
\label{Q2T2qAeven}
\end{figure}

If $q = 2t$ is even and $t>0$, then a presentation of the involutory quandle from Figure~\ref{torus diagrams with axis} is
$$Q_2(T_{2,q} \cup A) = \langle a, b, c \mid c^{ab} = c,\ a^{(ba)^{t-1} b c} = a,\ b^{(ab)^t c} = b\rangle_2.$$
As in the odd case, the Cayley graph is easier to produce using the following equivalent presentation
$$Q_2(T_{2,q} \cup A) = \langle a, b, c \mid \{c^{(ab)^i} = c \}_{i=1}^{2t-1},\ a^{(ba)^{t-1} c} = a,\ b^{(ab)^{t-1} c} = b,\ a^{cac}=a,\ b^{cbc}=b \rangle_2.$$
Applying the diagramming method to this presentation we find the Cayley graph shown in Figure~\ref{Q2T2qAeven}. The edges labeled $a(b)$ are $a$ when $t$ is even and $b$ when $t$ is odd (and the reverse for the edges labeled $b(a)$). As in the odd case, the order of the involutory quandle is $2k+2$.

We conclude by determining the automorphism groups for  $Q_2(T_{2,q} \cup A)$.

\begin{proposition} For $q \ge 2$, ${\rm Aut}(Q_2(T_{2,q} \cup A)) \cong \Z_2 \times {\rm Aut}(Q_2(L_{1/2q})) \cong \Z_2 \times \left(\mathbb Z_{2q} \rtimes \mathbb Z_{2q}^*\right)$.
\end{proposition}

{\bf Proof.} We will consider the case when $q = 2t+1$ is odd; the case when $q$ is even is similar.  From the presentation in  Section~\ref{2bridge}, notice that there is a natural inclusion $\eta : Q_2(L_{1/2q}) \rightarrow Q_2(T_{2,q} \cup A)$ determined by $\eta (a) =a$ and $\eta (b)=b$. Using the relation $x^{c(ab)^j}=x^{(ab)^jc}$ discussed in the derivation of second presentation of  $Q_2(T_{2,q} \cup A)$ we see that
$$\eta (a^{(ba)^{2t}b} )=a^{(ba)^{2t}b} \stackrel{R_3}{=} b^{c(ba)^{t}b} \stackrel{R_5}{=} b^{cb(ba)^{t}b} = b^{c(ab)^{t}}=b^{(ab)^{t}c} \stackrel{R_2}{=} a^{cc} =a = \eta(a).$$
In a similar manner, $\eta$ preserves the second relation $b^{(ab)^{2t}a}=b$ of $Q_2(L_{1/2q})$.  So, $\eta$ is a quandle homomorphism.  Clearly, $\eta$ is onto the set $H = \{x_1,\dots, x_{2t+1}, y_1, \dots , y_{2t+1}\}$ which has the same order, $2q$, as $Q_2(L_{1/2q})$. Thus, $\eta: Q_2(L_{1/2q}) \rightarrow H$ is a quandle isomorphism onto the subquandle $H$ of $Q_2(T_{2,q} \cup A)$.  Notice that $H$ is the subquandle generated by $\{a,b\}$ and that the set $S = \{z_1,z_2\} = Q_2(T_{2,q} \cup A) \setminus H$ is also a subquandle.

For any $f \in {\rm Aut}(Q_2(L_{1/2q}))$, we define $\varphi_f \in {\rm Aut}(Q_2(T_{2,q} \cup A))$ by $\varphi_f(a) = \eta \circ f(a)$, $\varphi_f(b) = \eta \circ f(b)$ and $\varphi_f(c) = c$ (from now on, we will abuse notation and use $f$ in place of $\eta \circ f$, as long as there is no confusion).  We also define $\psi \in {\rm Aut}(Q_2(T_{2,q} \cup A))$ by $\psi(a) = a$, $\psi(b) = b$ and $\psi(c) = c^a$.  We will show that ${\rm Aut}(Q_2(T_{2,q} \cup A)) = \langle \psi, \varphi_f \rangle \cong \Z_2 \times {\rm Aut}(Q_2(L_{1/2q}))$.

It is clear that $\langle \varphi_f \rangle \cong {\rm Aut}(Q_2(L_{1/2q}))$.  Also $\psi^2 = id$, since $\psi^2(c) = \psi(c^a) = \psi(c)^{\psi(a)} = c^{aa} = c$, so $\langle \psi \rangle \cong \Z_2$.  Finally, we compare $\phi_f \circ \psi$ and $\psi \circ \phi_f$.  Note that since $f(a)$ (respectively $f(b)$) involves only the generators $a$ and $b$, $\psi(f(a)) = f(a)$ (respectively, $\psi(f(b)) = f(b)$).

\begin{align*}
\phi_f \circ \psi(a) &= f(a), & \psi\circ\phi_f(a) &= \psi(f(a)) = f(a),\\
\phi_f\circ \psi(b) &= f(b), & \psi\circ\phi_f(b) &= \psi(f(b)) = f(b), \\
\phi_f\circ \psi(c) &= \phi_f(c^a) = c^{f(a)}, & \psi\circ\phi_f(c) &= \psi(c) = c^a.
\end{align*}

$f(a) = a^w$ or $b^w$, where $w$ is a word in $a$ and $b$.  Without loss of generality, suppose $f(a) = a^w$.  Then $c^{f(a)} = c^{a^w} = c^{\bar{w}aw}$.  Since the word $\bar{w}aw$ has odd length, $c^{\bar{w}aw} = c^a$.  Hence $\phi_f\circ \psi = \psi \circ \phi_f$ for every $\phi_f$.  This implies that $\langle \psi, \phi_f \rangle \cong \langle \psi \rangle \times \langle \phi_f \rangle \cong \Z_2 \times {\rm Aut}(Q_2(L_{1/2q}))$.

Finally, suppose $\alpha \in {\rm Aut}(Q_2(T_{2,q} \cup A))$.  Then $\alpha$ fixes $H$ and $S$ setwise, since $H$ is the only subquandle of order $2q$.  The restriction of $\alpha$ to $H$ gives an automorphism of ${\rm Aut}(Q_2(L_{1/2q}))$; let $f = \alpha\vert_H$.  Since $\alpha$ also fixes $S$, $\alpha(c) = c$ or $c^a$.  If $\alpha(c) = c$, then $\alpha = \phi_f$; on the other hand, if $\alpha(c) = c^a$, then $\alpha = \psi\circ \phi_f$.  So every automorphism is in $\langle \psi, \phi_f\rangle$.  Hence ${\rm Aut}(Q_2(T_{2,q} \cup A)) = \langle \psi, \varphi_f  \rangle \cong \Z_2 \times {\rm Aut}(Q_2(L_{1/2q}))$. \hfill $\Box$

\begin{corollary}
For $q > 2$, ${\rm Inn}(Q_2(T_{2,q} \cup A)) \cong D_{2q/\gcd(2,q)}$ and ${\rm Trans}(Q_2(T_{2,q} \cup A)) \cong D_q$.
\end{corollary}

\begin{figure}[htbp]
\begin{center}
\begin{tabular}{cc}
\includegraphics*[height=2.25in]{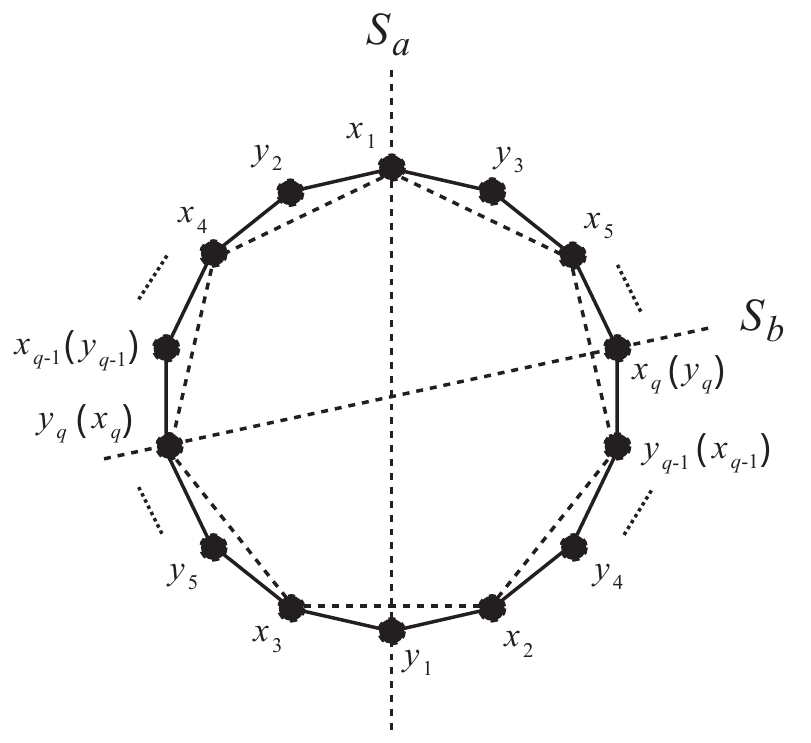}&
\includegraphics*[height=2.25in]{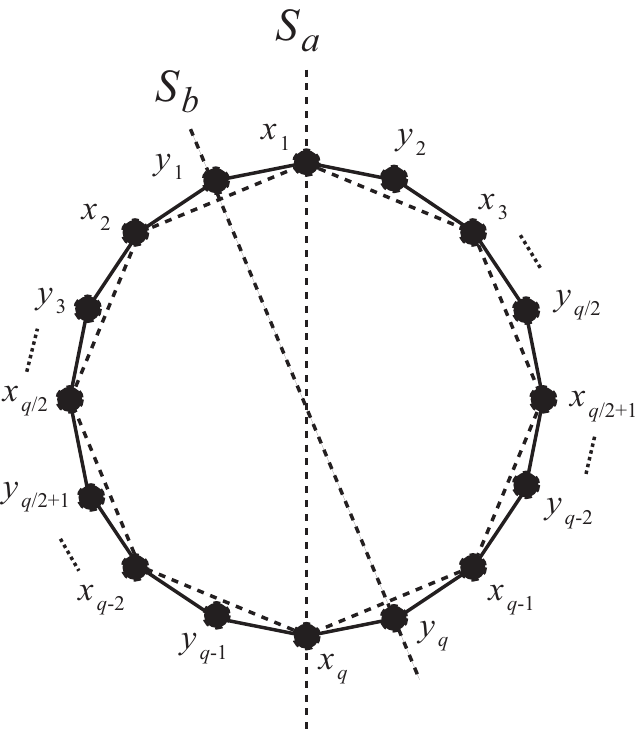}\\
\end{tabular}
\end{center}
\caption{Realizing $S_a$ and $S_b$ as symmetries for $q$ odd (on left) and even (on right).}
\label{polygons}
\end{figure}

{\bf Proof.} We will first consider the case when $q$ is odd.  The inner automorphism group is generated by the symmetries $S_a, S_b, S_c$.  We consider the Cayley graph from Figure \ref{Q2T2qJodd}, and use the vertex labelings from that diagram.  Then we can describe the action of each symmetry as a permutation on the set of vertices of the Cayley graph.
\begin{align*}
S_a &= (x_1)(x_2 x_3)(x_4 x_5) \cdots (x_{q-1} x_q) \cdot (y_1)(y_2y_3)(y_4y_5)\cdots (y_{q-1} y_q) \cdot (z_1 z_2), \\
S_b &= (x_1x_2)(x_3x_4) \cdots (x_{q-2} x_{q-1})(x_q) \cdot (y_1y_2)(y_3y_4) \cdots (y_{q-2} y_{q-1})(y_q)\cdot (z_1 z_2), \\
S_c &= (x_1y_1)(x_2y_2) \cdots (x_qy_q) \cdot (z_1)(z_2).
\end{align*}
These permutations correspond to symmetries of a regular $2q$-gon.  We embed the vertices of the $2q$-gon in the $xy$-plane, centered at the origin, and label them with $x_i$ and $y_i$ as shown on the left in Figure~\ref{polygons} (if $q \cong 1 \pmod{4}$, then $x_{q-1}$ and $y_q$ are on the left side; if $q \cong 3 \pmod{4}$, then $x_{q-1}$ and $y_q$ are on the right).  The vertices $z_1$ and $z_2$ are embedded on the $z$-axis above and below the polygon.  Then the actions of $S_a$ and $S_b$ are $180^\circ$ rotations around the lines through $\{x_1, y_1\}$ and $\{x_q, y_q\}$, respectively, and the action of $S_c$ is the $180^\circ$ rotation about the $z$-axis.  $S_a$ and $S_b$ generate the symmetries of a $q$-gon (the dotted polygon on the left in Figure \ref{polygons}).  Combined with the half-turn rotation of $S_c$ (since $q$ is odd), this generates the group of symmetries of the $2q$-gon.  So the inner automorphism group is isomorphic to $D_{2q}$.

Similarly, when $q$ is even, we consider the Cayley graph from Figure \ref{Q2T2qAeven}.  Once again, we describe $S_a$, $S_b$ and $S_c$ as permutations of the vertices in the Cayley graph.
\begin{align*}
S_a &= (x_1)(x_2 x_3)(x_4 x_5) \cdots (x_{q-2} x_{q-1})(x_q) \cdot (y_1y_2)(y_3y_4) \cdots (y_{q-1} y_{q}) \cdot (z_1 z_2), \\
S_b &= (x_1x_2)(x_3x_4) \cdots (x_{q-1} x_q) \cdot (y_1)(y_2y_3)(y_4y_5)\cdots (y_{q-2} y_{q-1})(y_q)\cdot (z_1 z_2), \\
S_c &= (x_1x_q)(x_2x_{q-1})\cdots (x_{q/2}x_{q/2+1}) \cdot (y_1y_q)(y_2y_{q-1})\cdots (y_{q/2}y_{q/2+1}) \cdot (z_1)(z_2).
\end{align*}
As before, we embed the vertices as the vertices of a $2q$-gon embedded in the $xy$-plane and centered at the origin, as shown on the right in Figure \ref{polygons} (once again, $z_1$ and $z_2$ are embedded on the $z$-axis; $x_{q/2}$ is on the left side if $q/2$ is even, and on the right if $q/2$ is odd).  As in the odd case, $S_a$ and $S_b$ generate the symmetries of the $q$-gon (the dotted polygon on the right in Figure \ref{polygons}); however, in this case the half-turn rotation $S_c$ is already in this group of symmetries.  So when $q$ is even, the inner automorphism group is just $D_q$. 

Since we are in an involutory quandle, each of $S_a$, $S_b$ and $S_c$ have order 2, and the transvection group is generated by their products $S_aS_b$, $S_aS_c$, and $S_bS_c$ (and their inverses).  Regardless of whether $q$ is odd or even, $S_aS_b$ is a rotation which generates the rotation subgroup of $D_q$, and $S_aS_c$ and $S_bS_c$ are reflections that have the same angle between their axes as $S_a$ and $S_b$.  So in both cases, these motions generate the symmetries of a $q$-gon, and the transvection group is isomorphic to $D_q$. \hfill $\Box$

We summarize our results on automorphism groups of torus links with axes in Table \ref{trefoil with axis automorphism groups}.

\begin{table}[h]
\begin{center}
\renewcommand{\arraystretch}{1.5}
\begin{tabular}{|c||c|c|c|c|}
\hline
$Q_n$ & $|Q_n|$ & \textrm{Aut} & \textrm{Inn} & \textrm{Trans} \\ \hline
$Q_2(T_{2,q} \cup A)$ & $2+2|q|$ &$\Z_2 \times \left(\mathbb Z_{2q} \rtimes \mathbb Z_{2q}^*\right) $ & $ D_{2q/\gcd(2,q)}$ &$D_q$ \\
\hline
$Q_2(T_{2,3} \cup B)$ & 18 & $\Z_2 \times \Z_2 \times S_4$ & $S_4$ & $S_4$ \\
\hline 
\end{tabular}
\end{center}
\caption{The automorphism groups of $Q_2(T_{2,3} \cup A)$ and $Q_2(T_{2,3} \cup B)$.}
\label{trefoil with axis automorphism groups}
\end{table}

{\bf Remark:} Aside from the 2-bridge links $L_{p/q}$, for which the involutory quandle depends only on $q$, all of the links in this paper are distinguished by the finite $n$-quandles described.  This can be easily seen by comparing the data in Tables \ref{automorphism groups of 2-bridge links}, \ref{pqnforfiniteQ}, \ref{torus quandle automorphism groups}, and \ref{trefoil with axis automorphism groups}.

 \end{document}